\documentclass[11pt,a4paper]{article}

\usepackage{latexsym}
\usepackage{blindtext}

\PassOptionsToPackage{force}{filehook}

\usepackage[utf8]{inputenc}
\usepackage[english]{babel}
\usepackage{arabtex}
\usepackage{utf8}
\usepackage{graphicx}
\usepackage{graphics,subfigure}
\usepackage{fancyhdr}
\usepackage{lmodern}
\usepackage{color}
\usepackage{transparent}
\usepackage{float}
\usepackage{pythontex}
\usepackage{tikz-cd}
\usepackage{listings}
\usepackage[linktocpage=true,citebordercolor=green, urlcolor=black]{hyperref}
\usepackage{tabularx}
\usepackage[title,titletoc]{appendix}
\usepackage{xcolor}
\usepackage[mathscr]{eucal}
\usepackage{calrsfs}
\usepackage{cite}
\usepackage{comment}

\usetikzlibrary{shapes.misc}
\usepackage{pgf,tikz}
\usepackage{caption}
\usetikzlibrary{decorations.markings}
\usetikzlibrary{arrows}
\usetikzlibrary{patterns}
\usepackage{tikz-3dplot}
\usetikzlibrary{calc}
\usepackage{amsfonts}
\usepackage{amsmath,amssymb,amsthm}
\usepackage{mathtools}
\usepackage{relsize}
\usepackage{titlesec}
\usepackage{arydshln}

\usepackage{stmaryrd}
\usepackage{trimclip}

\usepackage{enumitem}

\usepackage{wasysym}

\usepackage{cleveref}

\makeatletter
\DeclareRobustCommand{\shortto}{%
	\mathrel{\mathpalette\short@to\relax}%
}

\newcommand{\short@to}[2]{%
	\mkern2mu
	\clipbox{{.5\width} 0 0 0}{$\m@th#1\vphantom{+}{\shortrightarrow}$}%
}
\makeatother

\theoremstyle{definition}
\newtheorem{definition}{Definition}[section]

\newtheorem{remark}[definition]{Remark}

\newtheorem{lemma1}[definition]{Lemma}
\newtheorem{theorem1}[definition]{Theorem}

\newtheorem*{proposition*}{Proposition}

\titleformat{\paragraph}
{\normalfont\normalsize\bfseries}{\theparagraph}{1em}{}
\titlespacing*{\paragraph}
{0pt}{3.25ex plus 1ex minus .2ex}{1.5ex plus .2ex}

\usepackage[margin=2.5cm]{geometry}

\usepackage{mathptmx}

\AtEndDocument{%
	\par
	\medskip
	\begin{tabular}{@{}l@{}}%
		\textsc{Cengiz Aydin}\\
		\textsc{Institut für Mathematik, Universität Heidelberg, Germany}\\
		\textit{E-mail address}: \texttt{\href{mailto:cengiz.aydin@hotmail.de}{cengiz.aydin@hotmail.de}}
\end{tabular}}

\usepackage{authblk}
\title{Comet-type periodic motions and their out-of-plane bifurcations in the Earth--Moon CR3BP:\ a computational symplectic analysis}
\author{Cengiz Aydin}
\begin{document}
	
\setcounter{page}{1}
\pagenumbering{arabic}
	
\maketitle
	
\begin{abstract}
Comet-type periodic orbits of the circular restricted three-body problem (CR3BP) are periodic solutions that are generated from very large retrograde and direct circular Keplerian motions around the common center of mass of the primaries.\ In this paper we first provide an analytical proof of the existence of the comet-type periodic orbits by using the classical Poincaré continuation method.\ Within this analytical approach, we also determine the Conley--Zehnder index, defined as a Maslov index using a crossing form.\ Then, by applying a standard corrector--predictor technique, we explore numerically the two families of comet orbits within the Earth--Moon CR3BP.\ We compute their stability indices, identify vertical self-resonant bifurcations of higher order periods (of multiplicity from integer multiples up to six), investigate the vertically bifurcated spatial periodic solutions, and discuss their orbital characteristics.\ We also describe the orbits that are in resonance with the Earth and the Moon.\ We illustrate our main results in the form of bifurcation graphs, based on symplectic invariants, that provide a topological overview of the connections of the bifurcated branches, including bridge families.\
\end{abstract}

\begin{center}
\begin{tabular}{ll}
	\textbf{Keywords} & $\quad$ comet orbits $\cdot$ Earth--Moon CR3BP $\cdot$ vertical self-resonant bifurcations\\
	& $\quad$ bifurcation graph $\cdot$ Conley--Zehnder index $\cdot$ spatial resonant orbits\\
	\textbf{MSC 2020} & $\quad$ 70H12 $\cdot$ 70F07 $\cdot$ 70G45
\end{tabular}
\end{center}
		
\tableofcontents
		
\section{Introduction}

One of the most promising problems in celestial mechanics is the search for periodic solutions of the Newtonian $N$-body problem, which describes the dynamics of $N$ bodies of mass $m_i > 0$ ($i = 1,...,N$) moving and interacting according to Newton's law of gravitation.\ The equations of motion are given by
$$ \ddot{r}_i = - \sum_{j=1,j \neq i}^{N} \frac{m_j (r_i - r_j)}{r_{ij}^3},\quad 1 \leq i \leq N, $$
where the unit of time is chosen such that the gravitational constant is one, $r_i \in \mathbb{R}^3$ indicates the position of the $i$-body, and $r_{ij} = |r_i - r_j|$ denotes the Euclidean distance between the $i$-body and the $j$-body.\ While the 2-body problem has been completely understood \cite{roy}, the $N$-body problem for $N \geq 3$ still remains without a general analytical solution \cite{siegel_moser}.\ A particular and simpler case of interest is the restricted problem which describes conveniently the dynamics of small bodies of the solar system.\ One assumes that one of the $N$ bodies is significantly smaller than the other $N-1$ bodies (called primaries) and has therefore a negligible effect on their motion.\ With this assumption, the motion of the primaries becomes~a $(N-1)$-body problem to whose gravitational field the motion of the infinitesimal body is subject.\

Among the $N$-body systems, the 3-body problem is one of the most celebrated scenarios in each of astronomy, physics and mathematics.\ The circular restricted 3-body problem (CR3BP), in which the two primaries move in a circular fashion around their common center of mass, is used as a first approximation for the motion of a spacecraft in a Sun--planet or planet--moon environment.\ The periodic solutions of the CR3BP provide a paramount framework for the structure of all dynamically possible orbits and help to gain information on the stability and behavior of a given trajectory.\ In particular, the renewed popularity of Lunar explorations, spurred by NASA's Cislunar mission Artemis~\cite{williams}, along with the Lunar Gateway \cite{johnson}, and Chang'e missions \cite{liu}, demonstrate the crucial role of periodic orbits for operating orbits within space mission profiles and spacecraft trajectory designs.\ Therefore, it becomes of vital importance to catalog as many periodic orbits as possible.\

Although different kind of orbits of the Earth--Moon CR3BP have been widely studied \cite{broucke, doedel, zimovan, franz_russell}, there are still various types of periodic solutions that are sparsely explored in the literature, including the comet-type periodic orbits and their out-of-plane bifurcations.\ The comet-type periodic orbits are generated from very large retrograde and direct circular Keplerian motions around the barycenter of the primaries.\ Analytically, a first proof of their existence was given by Moulton \cite{moulton} for the R3BP.\ Later, Meyer provided an existence proof for (restricted) $N$-body problems \cite{meyer_1,meyer_2,meyer_3}, which is extended in \cite{llibre_roberto,llibre_stoica} for a wider class of restricted $N$-body problems.\ More recently, \cite{llibre_pasca_valls} proved the existence of comet orbits within a special CR4BP, and \cite{cors_garrido} analyzed them in an elliptic spatial restricted $N$-body problem.\ The technique of the existence proof is based on Poincaré's small parameter method, the Poincaré continuation method \cite{meyer_hall_offin}:\ for systems depending on a small parameter, a periodic solution of the unperturbed system exists in the perturbed system if its Floquet multipliers are different from the unity, as follows by the implicit function theorem.\

In classical works on numerical explorations, the comet orbits and their planar linear stability were investigated by Strömgren and his colleagues within the Copenhagen problem \cite{szebehely} (a special case of the CR3BP when the two primaries have equal mass), by Broucke \cite{broucke} (for the Earth--Moon CR3BP), and by Bruno \cite{bruno_1,bruno_2} (for the Sun--Jupiter system).\ The corresponding orbit family notations are:\
\begin{center}
\begin{tabular}{lccc}
	& CR3BP & direct comet orbits & retrograde comet orbits\\
	Strömgren & Copenhagen problem & $l$ & $m$\\
	Broucke & Earth--Moon system & $E_1$ & $F$\\
	Bruno & Sun--Jupiter system & $ID-1$ & $IR-$
\end{tabular}
\end{center}

The out-of-plane bifurcations from the comet orbits were analyzed in \cite{lara_scheeres} for a gravitational model on the dynamics around the asteroid (433) Eros, and in \cite{papadakis} for a particular case of the CR4BP.\ For various CR3BP systems, \cite{doedel} describes only partially some vertically bifurcated branches from retrograde comet orbits (labeled as family $C1$).\ Since we are not aware of any additional articles on the comet orbits in the literature, the present work aims at studying the retrograde and direct comet periodic orbits and their vertical multiple-cover bifurcations, with a particular focus on the Earth--Moon environment.\

In recent years, symplectic invariants (most significant is the Conley--Zehnder index) have been successfully applied to the numerical bifurcation analysis, leading to a topological overview of the global network structure and systematic classification of periodic orbit families \cite{aydin_cz,moreno_aydin,aydin_batkhin,aydin_dro,joung_otto,lee,joung_koh}.\ This work now explores the combination of these symplectic methods with classical techniques to investigate the comet orbit families of the Earth--Moon system.\ Our results show that in almost all cases the bifurcated branches form bridge families.\ By a \textit{bridge family} we understand a periodic orbit family between two different branch points, as defined in \cite{llibre}.\ Therefore, bridge families are \textit{closed families}, existing only on a finite interval of the energy values.\ Throughout the paper, we denote the family of retrograde and direct comet orbits by $\kappa_-$ and $\kappa_+$, respectively.\ Our vertical multiple-cover bifurcation results are:\
\begin{itemize}[noitemsep]
	\item[a)] \textit{Single-turn and period-doubling bifurcations from the retrograde comet orbits (family $\kappa_-$).}\\
	The solutions that are generated from the single-turn bifurcation from the family $\kappa_-$ correspond to the $L_1$ halo orbits, and those that are computed from the period-doubling bifurcation from the family $\kappa_-$ agree with the families of the $L_2$ and $L_3$ vertical Lyapunov orbits.\ We mention that these results are also described in \cite{doedel}, but are derived from the solutions around the libration points.\
	\item[b)] \textit{Period-tripling bifurcation from retrograde comet orbits (family $\kappa_-$).}\\
	The period-tripling bifurcation from the $\kappa_-$ orbits gives rise to two bifurcated branches.\ While one branch forms a bridge family to a periodic solution which belongs to the family of retrograde solutions around the Earth (family $A_1$ in Broucke's notation), the other branch forms a bridge to the single-turn bifurcation from the direct comet orbits (family $\kappa_+$).\
	\item[c)] \textit{Bifurcations from $k$-covered retrograde comet orbits (family $\kappa_-$), for $k=4,5,6$.}\\
	Systematically, for each $k=4,5,6$, all the bifurcated branches from the $k$-covering of the $\kappa_-$ orbits form bridge families to the $(k-2)$-covering of $\kappa_+$ orbits.\
	\item[d)] \textit{Single-turn bifurcation from the direct comet orbits (family $\kappa_+$)}.\\
	From the single-turn bifurcation of the $\kappa_+$ orbits, we generate a branch whose members along the continuation collide with the Earth surface.\ Most of the period the solutions are far from the Earth with a large excursion of hook-shaped form in the vertical direction, so that the near-collision phase is very brief.\ This kind of orbit is known within the CR3BP around the larger primary for small mass values, as bifurcation from vertical collision orbits in the rotating Kepler problem \cite{belbruno_1,belbruno_2}.\
\end{itemize}

\textit{Structure of the paper.}\ In Section \ref{sec2} we recall the CR3BP dynamics.\ The following Section \ref{sec3} contains the construction of the Conley--Zehnder index and associated symplectic invariants.\ In Section \ref{sec4} we provide an analytical approach to the comet orbits, where we first discuss an existence proof based on the Poincaré continuation method, and then determine their Conley--Zehnder indices.\ In Section \ref{sec5} we numerically explore the comet orbit families and present the outcomes.\ We compute stability indices, identify vertical self-resonant orbits up to the multiplicity six, analyze the vertically bifurcated branches and discuss their orbital characteristics.\ The main results, that are summarized briefly above, are illustrated in the form of bifurcation graphs.\ We also discuss the orbits that are in resonance with the Earth and the Moon.\ In Section~\ref{sec6} we summarize our work and conclude.\ Data base are provided in the Appendix.\

\section{CR3BP dynamics}
\label{sec2}

\subsection{Equations of motion}

\begin{figure}[t!]
	\centering
	\includegraphics[scale=1]{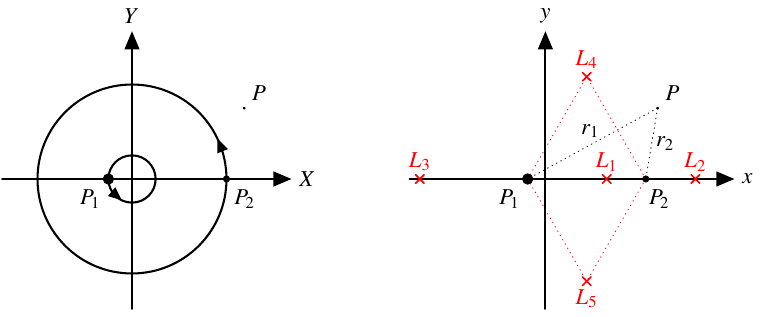}
	\caption{CR3BP (from \cite{aydin_dro}) in an inertial reference frame $(X,Y)$ and in a rotating reference frame $(x,y)$, with primaries $P_1$ and $P_2$, infinitesimal body $P$, libration points $L_i$ ($i=1,2,3,4,5$), and distances $r_1$ and~$r_2$ from $P$ to the primaries $P_1$ and $P_2$, respectively.}
	\label{figure_1_cr3bp}
\end{figure}

The CR3BP refers to the dynamics of a body $P$ of infinitesimal mass that is subject to the gravitational field of two primary bodies, a larger primary $P_1$ and a smaller primary $P_2$ with masses $m_1 \geq m_2$, that move along circles in the same plane with constant angular velocity around their common center of mass (see Figure \ref{figure_1_cr3bp}).\ The infinitesimal body $P$ has a negligible effect on the motion of the primaries.\ Taking a coordinate system that rotates with an angular velocity of the orbital angular rate of the primaries, and choosing suitable units, the primaries $P_1$ and $P_2$ can be assumed to have respective masses $1 - \mu$ and~$\mu$ ($\mu \in [0,\frac{1}{2}]$), to rest at the respective points $(-\mu,0,0)$ and $(1-\mu,0,0)$ (see Figure \ref{figure_1_cr3bp}).\ Moreover, the primaries $P_1$ and $P_2$ complete one inertial revolution in $2 \pi$ time units.\ Within these assumptions, the motion of the infinitesimal body is governed by the second order differential equations \cite[Chapter 10]{szebehely}:\
\begin{align}
	\ddot{x} &= 2 \dot{y} + x - (1 - \mu) \dfrac{x + \mu}{r_1^3} - \mu \dfrac{x - 1 + \mu}{r_2^3}, \nonumber \\
	\ddot{y} &= - 2 \dot{x} + y - \left( \dfrac{1 - \mu}{r_1^3} + \dfrac{\mu}{r_2^3} \right) y, \label{ham_equation} \\
	\ddot{z} &= - \left( \dfrac{1-\mu}{r_1^3} + \dfrac{\mu}{r_2^3} \right) z, \nonumber
\end{align}
where $r_1 = \left( (x+\mu)^2 + y^2 + z^2 \right)^{\frac{1}{2}}$ and $r_2 = \left( (x - 1 + \mu)^2 + y^2 + z^2 \right)^{\frac{1}{2}}$ indicate the distances from the infinitesimal body to the primaries.\ The CR3BP framework has five points of equilibrium, also called libration points or Lagrange points \cite[Chapter 4]{szebehely}:\ Three collinear points $L_1$ (located within the primaries), $L_2$ and $L_3$ (outside the interval joining the primaries), and two equilateral triangular points $L_4$ and $L_5$, as shown in Figure \ref{figure_1_cr3bp}.\

Defining kinetic moments by
\begin{align}\label{kinetic_moments}
	p_x = \dot{x} - y,\quad p_y = \dot{y} + x, \quad p_z = \dot{z},
\end{align}
the system can be written in Hamiltonian form with the corresponding Hamiltonian function
\begin{align} \label{hamiltonian}
	H(x,y,z,p_x,p_y,p_z) = \frac{1}{2} \left( p_x^2 + p_y^2 + p_z^2 \right) - \frac{1-\mu}{r_1} - \frac{\mu}{r_2} + p_x y - p_y x,
\end{align}
which is a first integral of the system.\ Thus, periodic solutions come in a smooth one-parameter family, parameterized by the energy.\ An equivalent first integral is the Jacobi constant $C = -2H$.\ The phase space of the system is 6-dimensional and endowed with the standard symplectic form $\omega_0 = \sum dp_k \wedge dk$ ($k=x,y,z$).\ The Hamiltonian vector field $X_H$, uniquely determined by $dH(\cdot) = \omega_0( \cdot, X_H )$, generates the Hamiltonian flow $\phi^t$ by the equations
\begin{align}\label{equations_2}
	\frac{d}{dt} \phi^t = X_H (\phi^t),\quad \phi^0 = \text{id},
\end{align}
which is equivalent to \eqref{ham_equation}.\ In the linear algebra expression of $\omega_0$,
$$ \omega_0 (\cdot, \cdot) = \langle J \cdot, \cdot \rangle,\quad \text{with } J = \begin{pmatrix}
	0 & I_3\\
	-I_3 & 0
\end{pmatrix}, $$
where the brackets denote the Euclidean inner product, we can represent the Hamiltonian vector field by
$$X_H = J \nabla H,$$
where the gradient, as usual, is defined with respect to the Euclidean inner product.\ With that, the equations (\ref{ham_equation}) and (\ref{equations_2}) can be equivalently expressed by
\begin{align}\label{equation_J}
	\dot{ \gamma } = J \nabla H ( \gamma ),\quad \gamma \equiv (x,y,z,p_x,p_y,p_z).
\end{align}

\subsection{Monodromy matrix and stability index}

Let $\gamma \equiv (x,y,z,p_x,p_y,p_z)$ be a $T$-periodic solution of (\ref{equation_J}).\ The linearized flow, which is represented by
\begin{align*}
	d \phi^t =: V(t) = (v_{ij}) \in \mathbb{R}^{6 \times 6},
\end{align*}
solves the \textit{variational equation},
\begin{align}\label{var_eq}
	\dot{V}(t) = J D^2 H \left( \gamma (t) \right) V(t),\quad V(0) = \text{id},
\end{align}
where $D^2 H $ is the Hessian of the Hamiltonian.\ The \textit{monodromy matrix}~$M$ of~$\gamma$ corresponds to the full period linearization, i.e.,\ $M = V(T)$.\ Since the monodromy matrix is symplectic (it satisfies $M^T J M = J$), or it is conjugate to a symplectic matrix by (\ref{kinetic_moments}), the eigenvalues occur in reciprocal pairs and are of the form
$$\{1,1,\lambda_1,\lambda_1^{-1},\lambda_2,\lambda_2^{-1}\},$$
where 1 appears trivially twice since the energy is a first integral of the system.\ The non-trivial eigenvalues are called \textit{Floquet multipliers}.\ The periodic orbit $\gamma$ is \textit{non-degenerate} when~1 is not among its Floquet mutlipliers.\ For each pair of Floquet multipliers the \textit{stability index} is defined as
$$ s_k = \frac{1}{2} \left( \lambda_k + \lambda_k^{-1} \right),\quad \text{for } k = 1,2.$$
They can be of one of the following types:
\begin{itemize}[noitemsep]
	\item \textit{Elliptic:}\ $-1 < s_k = \cos \theta < 1$.\ It is equivalent to $\lambda_k = e^{i \theta}$.
	\item \textit{Positive hyperbolic:}\ $s_k > 1$.\ It is equivalent to $\lambda_k \in \mathbb{R}_{>0} \setminus \{1\}$.
	\item \textit{Negative hyperbolic:}\ $s_k < -1$.\ It is equivalent to $\lambda_k \in \mathbb{R}_{<0} \setminus \{-1\}$.
	\item \textit{Complex instability:}\ Case of quadruples $\{ \lambda, \lambda^{-1}, \overline{\lambda}, \overline{\lambda}^{-1} \}$.
\end{itemize}
Stable periodic solutions are those that simultaneously satisfy $|s_1| < 1$ and $|s_2| < 1$.\ It is common to compute the stability indices via formulas that are given in the next lemma following \cite{bray_goudas}.\
\begin{lemma1}
	\textit{The stability indices are determined by}
	\begin{align}\label{stab_ind}
		s_1 = - \frac{1}{4} \left( \alpha + \sqrt{\beta} \right),\quad s_2 = - \frac{1}{4} \left( \alpha - \sqrt{\beta} \right),
	\end{align}
	\textit{where $\alpha = 2 - \textnormal{tr}(M)$, and $\beta  = 2 \textnormal{tr}(M^2) - \alpha^2 + 4$.}
\end{lemma1}
\begin{proof}
	The characteristic polynomial of $M$ is of the form
	\begin{align}\label{char}
		\chi_M(\lambda) = (\lambda - 1)^2 \left( \lambda^4 + \alpha \lambda^3 + \beta' \lambda^2 + \alpha \lambda + 1 \right),
	\end{align}
	where
	$$\alpha = - \lambda_1 - \lambda_1^{-1} - \lambda_2 - \lambda_2^{-1},\quad \beta' = \lambda_1 \lambda_2 + \lambda_1 \lambda_2^{-1} + \lambda_1^{-1} \lambda_2 + \lambda_1^{-1} \lambda_2^{-1} + 2.$$
	The coefficients $\alpha$ and $\beta'$ can also be written as
	$$ \alpha = 2 - \text{tr}(M),\quad \beta' = \frac{1}{2} \alpha^2 - \frac{1}{2} \text{tr} (M^2) + 1. $$
	By taking into account the stability indices, (\ref{char}) can be equivalently written as
	\begin{align*}
		\chi_M (\lambda) = (\lambda - 1)^2 \left( \lambda^4 + 2 (- s_1 - s_2) \lambda^3 + (4s_1 s_2 + 2) \lambda^2 + 2(-s_1 - s_2)\lambda + 1 \right).
	\end{align*}
	Therefore, the stability indices satisfy the equations 
	$$ -2(s_1 + s_2) = \alpha,\quad 4 s_1 s_2 + 2 = \beta',$$
	whose solutions correspond to (\ref{stab_ind}).\
\end{proof}

\subsection{Symmetries}

\textit{Symmetries} are the following involutions on the phase space that leave \eqref{hamiltonian} invariant:
\begin{align}
	\rho_x (x,y,z,p_x,p_y,p_z) = &(x,-y,-z,-p_x,p_y,p_z) \text{ ($\pi$-rotation around the $x$-axis)},\label{rho1}\\
	\rho_{xz} (x,y,z,p_x,p_y,p_z) = &(x,-y,z,-p_x,p_y,-p_z) \text{ (reflection at the $xz$-plane)},\label{rho2}\\
	\sigma (x,y,z,p_x,p_y,p_z) = &(x,y,-z,p_x,p_y,-p_z) \text{ \textcolor{white}{$-$}(reflection at the ecliptic $\{z=0\}$)}.\label{sigma}
\end{align}
\noindent
Notice that $\rho_x \circ \rho_{xz} = \rho_{xz} \circ \rho_{x} = \sigma$.\ The planar problem corresponds to the restriction of the system to the fixed point set Fix$(\sigma) = \{ (x,y,0,p_x,p_y,0) \}$, in which an inherent symmetry is characterized by $(x,y,0,p_x,p_y,0) \mapsto (x,-y,0,-p_x,p_y,0)$ (reflection at the $x$-axis).\ Moreover, while (\ref{sigma}) is symplectic, (\ref{rho1}) and (\ref{rho2}) are anti-symplectic, which means that their matrix representations,
\begin{align*}
	\Sigma_{X} &= \text{diag}\{ 1,-1,- 1,-1,1,1 \},\\
	\Sigma_{XZ} &= \text{diag}\{ 1,-1,1,-1,1,- 1 \},\\
	\Sigma_{\sigma} &= \text{diag}\{ 1,1,-1,1,1,-1 \},
\end{align*}
satisfy $\Sigma_{\sigma} J \Sigma_{\sigma} = J$, and $\Sigma J \Sigma = - J$, where $\Sigma = \Sigma_{X}$, or $\Sigma =\Sigma_{XZ}$.\ While the symplectic symmetries denote the time-preserving symmetries, the anti-symplectic symmetries denote the time-reversal ones in the Hamiltonian context~\cite{lamb_roberts}.\ Furthermore, the Hamiltonian flow~$\phi^t$ satisfies
\begin{align}\label{identity}
	\sigma \circ \phi^t = \phi^t \circ \sigma,
\end{align}
and for the anti-symplectic symmetries, it holds that
\begin{align}\label{identity_2}
	\rho_i \circ \phi^{-t} = \phi^t \circ \rho_i, \quad \text{with } i \in \{x,xz\}.
\end{align}
In particular, the periodic orbits that are invariant under anti-symplectic symmetries are divided into groups with different number of symmetries:
\begin{itemize}[noitemsep]
	\item \textit{Simple symmetric w.r.t.\ $\rho_{x}$ (or $\rho_{xz}$)}:\ The orbit starts perpendicular at the $x$-axis (or $xz$-plane) and intersects the $x$-axis (or $xz$-plane) perpendicularly at half period.\ 
	\item \textit{Doubly symmetric w.r.t.\ $\rho_x$-$\rho_{xz}$ (or $\rho_{xz}$-$\rho_x$)}:\ The orbit starts perpendicular at the $x$-axis (or $xz$-plane) and intersects the $xz$-plane (or $x$-axis) perpendicularly at quarter period.
\end{itemize}

The analysis of symmetric orbits can be simplified by using the symmetries \cite{bray_goudas}:\ It is sufficient to integrate the equations from one symmetric point to the next symmetric point to obtain the monodromy matrix, which we formulate in the next lemma.\ This in fact, offers economy of computing effort.\

\begin{lemma1}\label{lemma_mon}
	\textit{Let $\gamma$ be a simple or doubly symmetric $T$-periodic orbit.\ Depending on the type of symmetries, the monodromy matrix of $\gamma$ is determined as follows:}
	\begin{itemize}[noitemsep]
		\item \textit{a) If $\gamma$ is simple symmetric w.r.t. $\rho_{x}$, or $\rho_{xz}$, then correspondingly:}
		$$ M = \Sigma_{X} V^{-1}(T/2) \Sigma_{X} V(T/2),\quad \text{or} \quad M = \Sigma_{XZ} V^{-1}(T/2) \Sigma_{XZ} V(T/2).$$
		\item \textit{b) If $\gamma$ is doubly symmetric w.r.t. $\rho_{x}$-$\rho_{xz}$, or $\rho_{xz}$-$\rho_{x}$, then correspondingly:}
		$$ M = \left[ \Sigma_{X} V^{-1}(T/4) \Sigma_{XZ} V(T/4) \right]^2,\quad \text{or} \quad M = \left[ \Sigma_{XZ} V^{-1}(T/4) \Sigma_{X} V(T/4) \right]^2.$$
	\end{itemize}
\end{lemma1}
\begin{proof}
	The proof for the respective second case in a) and b) is similar to that in the respective first cases, therefore it is sufficient to discuss only the corresponding first case.\ Let $\gamma$ be a simple symmetric solutions w.r.t.\ $\rho_{x}$.\ Then, in view of the identity (\ref{identity_2}), the variational matrix satisfies
	\begin{align}\label{iden_1}
		V^{-1}(t) = \Sigma_{X} V(t) \Sigma_{X}.
	\end{align}
	The equation (\ref{iden_1}) for $t = - T/2$, together with $V(T/2) = V^{-1}(T/2) V(T)$, becomes
	$$ V^{-1}(T/2) V(T) = \Sigma_{X} V^{-1}(T/2) \Sigma_{X}, $$
	which implies the first case in a).\ Now let $\gamma$ be a doubly symmetric solutions w.r.t.\ $\rho_{x}$-$\rho_{xz}$.\ Since the orbit $\gamma$ is symmetric w.r.t.\ $\rho_{x}$ on half of its period, from the first case in a) we have that
	\begin{align}\label{mon}
		M = \Sigma_{X} V^{-1}(T/2) \Sigma_{X} V(T/2).
	\end{align}
	Taking into account that $\Sigma_{X} \Sigma_{XZ} = \Sigma_{XZ} \Sigma_{X} = \Sigma_{\sigma}$, the variational matrix satisfies, in view of (\ref{identity}), that
	\begin{align}\label{iden_2}
		V(t + T) = \Sigma_{\sigma} V(t) \Sigma_{\sigma} V(T).
	\end{align}
	The combination of (\ref{iden_1}) and (\ref{iden_2}) for $t = -T/4$  gives
	\begin{align*}
		V(T/2) = \Sigma_{\sigma} \Sigma_{X} V^{-1}(T/4) \Sigma_{X} \Sigma_{\sigma} V(T/4) = \Sigma_{XZ} V^{-1}(T/4) \Sigma_{XZ} V(T/4),
	\end{align*}
	which, together with its inverse, inserted into (\ref{mon}) proves the first case in b).\
\end{proof}

\begin{remark}
	In symplectic coordinates, the variational matrix $V(t)$ satisfies $V^T(t) J V(t) = J$, which means that in the formulas in Lemma \ref{lemma_mon}, the inverse can be determined by $V^{-1}(t) = J^T V^T(t) J$.\
\end{remark}

\subsection{Out-of-plane bifurcations}

By virtue of the symplectic symmetry $\sigma$ (\ref{sigma}), the planar periodic orbits have a planar and spatial stability index (also called \textit{horizontal} or \textit{vertical stability index}), that we denote correspondingly by $s_p$ and $s_v$.\ The case $s_v = 1$ signals an out-of-plane bifurcation with spatial orbits of the same period, while the case $s_v = -1$ corresponds to a period-doubling bifurcation.\ The \textit{vertical self-resonant} (VSR) orbits are those with Floquet multipliers equal to the root of the unity, i.e.,\ the vertical stability index satisfies
\begin{align}\label{vsr}
	s_v = \cos \left( 2 \pi \frac{d}{k} \right),
\end{align}
where $d,k \in \mathbb{Z}$ $(k \neq 1,2)$ with $d<k$ and $d/k$ irreducible fraction.\ Such a critical orbit is also called $d\text{:}k$ \textit{self-resonant orbit}.\ As identified in \cite{robin_markellos}, the $k$-covering with period $kT$ generates bifurcation of exactly two spatial periodic orbit families whose symmetry properties depend on the multiplicity $k$:
\begin{itemize}[noitemsep]
	\item \textit{Odd multiplicity $k$}.\ The members of the two bifurcated branches are simple symmetric; while the orbits of one family are symmetric w.r.t.\ $\rho_{x}$, the others are symmetric w.r.t.\ $\rho_{xz}$.\ Each branch exists twice by using the $\sigma$-symmetry.\
	\item \textit{Even multiplicity $k$}.\ The orbits of each of the bifurcated branches are doubly symmetric w.r.t.\ $\rho_{x}$ and $\rho_{xz}$.\
\end{itemize}

\section{Conley--Zehnder index}
\label{sec3}

\subsection{Transverse Conley--Zehnder index of a Hamiltonian orbit}

An important invariant associated to Hamiltonian periodic solutions is the so-called Conley--Zehnder index \cite{conley_zehnder}.\ Intuitively, it measures the winding number of the linearized flow along an orbit.\ Due to the time independence of the Hamiltonian framework we study, we consider the transverse Conley--Zehnder index of a periodic solution such as proposed in \cite{albers,moreno_aydin,joung_otto,lee,joung_koh}.\ For our purpose, we look attentively at the approach in \cite{albers,lee} which is based on the Robbin--Salamon definition of the Maslov index using a crossing form \cite{robbin_salamon}:\ Let $Sp(2n)$ denote the symplectic group, i.e.,\
$$Sp(2n) = \{ A \in \text{Mat}(2n,\mathbb{R}) \colon A^T J A = J \},\quad \text{with } J = \begin{pmatrix}
	0 & I_n\\
	-I_n & 0
\end{pmatrix}.$$
Let $\Psi \colon [0,T] \rightarrow Sp(2n)$ be a path of symplectic matrices with $\Psi(0) = \text{id}$.\ A point $t$ is called a \textbf{crossing} if $\text{ker} \left( \Psi(t) - \text{id} \right) \neq 0$.\ For a crossing $t$ the \textbf{crossing form} is defined as a quadratic form on the vector space $V_t = \text{ker} \left( \Psi(t) - \text{id} \right) $ by
$$ Q_t(v,v) = \omega \left( v, \dot{\Psi}(t) v \right),\quad v \in V_t, $$
where $\omega$ is a symplectic form on $V_t$.\ Since in symplectic coordinates, $\omega$ is represented by the matrix
\begin{align*}
	\Omega := \text{diag} \left\{ \begin{pmatrix}
		0 & 1\\
		-1 & 0
	\end{pmatrix}, ..., \begin{pmatrix}
		0 & 1\\
		-1 & 0
	\end{pmatrix} \right\},
\end{align*}
the crossing form can be expressed as a diagonal matrix by $Q_t = \Omega \dot{\Psi}(t)$.\ Furthermore, the signature of the quadratic form is defined as $\text{sign}(Q_t) = n_+ - n_-$, where $n_+$ and $n_-$ denote the counts of positive and negative entries, respectively.\ Assume that all crossings are isolated and non-degenerate, i.e., the crossing form $Q_t$ at the crossing $t$ is non-degenerate as a quadratic form.\ The \textbf{Maslov index} of $\Psi$ (also called \textbf{Robbin--Salamon index}) is then defined as
\begin{align}\label{maslov}
	\mu (\Psi) = \frac{1}{2} \text{sign} (Q_0) + \sum_{t \in (0,T) \text{ crossing}} \text{sign}(Q_t) + \frac{1}{2} \text{sign}(Q_T),
\end{align}
which has the following nice properties:
\begin{itemize}[noitemsep]
	\item[1)] (\textit{Homotopy}) Two paths of symplectic matrices are homotopic with end points fixed if and only if they have the same index.
	\item[2)] (\textit{Catenation}) If $\Psi \colon [a,b] \rightarrow Sp(2n)$ and $a < c < b$, then
	\begin{align*}
		\mu(\Psi) = \mu(\Psi | _{[a,c]}) + \mu(\Psi |_{[c,b]}).
	\end{align*}
	\item[3)] (\textit{Product}) Let $n = n' + n''$, and identify $Sp(2n') \times Sp(2n'')$ as a subgroup of $Sp(2n)$ in the obvious way.\ If $\Psi', \Psi''$ are paths in $Sp(2n')$, $Sp(2n'')$ respectively, then
	\begin{align}\label{index_product}
		\mu(\Psi' \oplus \Psi'') = \mu(\Psi') + \mu(\Psi'').
	\end{align}
\end{itemize}

\begin{figure}[t!]
	\centering
	\includegraphics[scale=1]{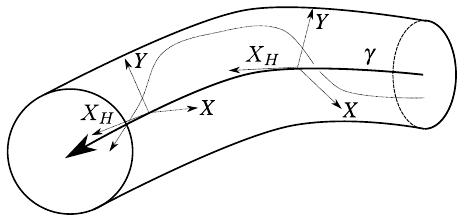}
	\caption{Conley--Zehnder index measures a twisting of the linearized flow along an orbit $\gamma$ with respect to a frame (reproduced from \cite{moreno_aydin} with minor modifications).}
	\label{figure_2_cz}
\end{figure}

Now let $H$ be a time-independent Hamiltonian on a phase space, say $\mathbb{R}^{2n}$, and let $\gamma$ be a $T$-periodic orbit of the Hamiltonian vector field $X_H$.\ In order to relate the Maslov index to the orbit $\gamma$, we need a path of symplectic matrices, which in fact, is generated by the linearized flow along the orbit $\gamma$.\ We first have to fix a symplectic frame, or a symplectic trivialization, along the orbit $\gamma$.\ This frame consists of a symplectic basis of the tangent space at each point of the orbit $\gamma$, with respect to which we measure the winding of the linearized flow, as illustrated in Figure \ref{figure_2_cz}.\ The existence of such a symplectic trivialization is guaranteed by general theory \cite[Chapter 2.6]{mcduff_salamon} within the assumption that the orbit $\gamma$ is the boundary of a spanning disk in~$\mathbb{R}^{2n}$.\

Since the preservation of the energy, we first consider at each point of the orbit $\gamma$ a 2-dimensional symplectic vector space determined by
\begin{align*}
	L := \langle Z, X_H \rangle,\quad \text{with } Z = \frac{1}{|| \nabla H ||^2} \nabla H.
\end{align*}
The normalization of the vector $Z$ ensures that $\omega(Z,X_H) = 1$.\ Therefore, $\{Z,X_H\}$ form a symplectic basis of the frame $L$.\ We next have to choose a symplectic basis $\{ X_1, ... , X_{n-1}, Y_1, ... , Y_{n-1}\}$ of the $(2n-2)$-dimensional symplectic complement of $L$, that is defined by
\begin{align}\label{sym_complement}
	L^{\omega} = \{ Y \colon \omega(X,Y) = 0 \text{ for } X \in L \}.
\end{align}
Such a symplectic basis can be achieved, for example, by the symplectic Gram--Schmidt process.\ In this way, on an energy level set $\Sigma = H^{-1}(c)$, $c \in \mathbb{R}$, at every point $p$ of the orbit $\gamma$ the symplectic basis of $L^{\omega}$ forms a $(2n-2)$-dimensional frame that is transverse to the Hamiltonian vector field $X_H$.\ This frame determines a symplectic trivialization by
$$ \Phi \colon \Sigma \times \mathbb{R}^{2n-2} \to L^{\omega},\quad (p,u_1,...,u_{n-1},v_1,...,v_{n-1}) \mapsto \sum_{i=1}^{n-1}u_i \cdot X_i(p) + \sum_{i=1}^{n-1}v_i \cdot Y_i(p).$$
Together with the projection map,
$$ pr_{L^{\omega}} \colon T\mathbb{R}^{2n} \to L^{\omega},\quad a Z + \sum_{i=1}^{n-1}u_i X_i + b X_H + \sum_{i=1}^{n-1}v_i Y_i \mapsto \sum_{i=1}^{n-1}u_i X_i + \sum_{i=1}^{n-1}v_i Y_i, $$
a path of symplectic matrices $\Psi_{\gamma} \colon [0,T] \to Sp(2n-2)$, that represents the transverse linearized flow, is now generated by
\begin{align}\label{path}
	\Psi_{\gamma}(t) := \left( \Phi(\gamma(t),\cdot) \right)^{-1} \circ pr_{L^{\omega}} \circ d \phi^t(\gamma(0)) \circ \Phi(\gamma(0),\cdot) ,
\end{align}
which is characterized by the diagram:\
$$
\begin{tikzcd}[column sep=10em]
	\mathbb{R}^{2n-2} \arrow[r,dashed, "\Psi_{\gamma}(t)"] \arrow[swap,d, "{\Phi(\gamma(0),\cdot)}"] & \mathbb{R}^{2n-2} \\
	L_{\gamma(0)}^{\omega} \arrow[r,"pr_{L^{\omega}} \circ d \phi^t(\gamma(0))"] & L_{\gamma(t)}^{\omega} \arrow[swap,u, "{\left( \Phi(\gamma(t),\cdot) \right)^{-1}}"]
\end{tikzcd}
$$
Then, the \textbf{transverse Conley--Zehnder index} of the orbit $\gamma$ is defined as the Maslov index of the path~$\Psi_{\gamma}$, which is invariant under continuous deformations of trivializations and that we denote by
$$ \mu_{CZ} (\gamma) := \mu(\Psi_{\gamma}). $$

\subsection{Planar solutions, the case of $Sp(2)$ and index iteration}

Let $\gamma \equiv (x,y,0,p_x,p_y,0)$ be a planar $T$-periodic solution of a spatial Hamiltonian system in which the planar framework occurs as the fixed point set of the symplectic symmetry $\sigma$ (\ref{sigma}).\ Then, along the planar orbit $\gamma$ the 4-dimensional transverse frame splits symplectically into 2-dimensional planar and spatial frames, which we denote correspondingly by $[X_p \text{ } Y_p]$ and $[X_s \text{ } Y_s]$.\ In this setting, a global trivialization can be chosen quite simply and explicitly:\ For the spatial frame, the out-of-plane directions,
\begin{align}\label{spatial_frame}
	X_s = \partial_{p_z},\quad Y_s = \partial_z,
\end{align}
form a symplectic basis since $\omega(X_s,Y_s)=1$.\ For the planar frame, it is convenient to work with quaternionic matrices \cite{moreno_aydin,joung_otto}, which, restricted to planar variables, are given by
$$ J = \begin{pmatrix}
	0 & 0 & 1 & 0\\
	0 & 0 & 0 & 1\\
	-1 & 0 & 0 & 0\\
	0 & -1 & 0 & 0
\end{pmatrix},\quad K_1 = \begin{pmatrix}
0 & 1 & 0 & 0\\
-1 & 0 & 0 & 0\\
0 & 0 & 0 & -1\\
0 & 0 & 1 & 0
\end{pmatrix},\quad K_2 =  \begin{pmatrix}
0 & 0 & 0 & -1\\
0 & 0 & 1 & 0\\
0 & -1 & 0 & 0\\
1 & 0 & 0 & 0
\end{pmatrix}. $$
 The vectors
 \begin{align}\label{planar_frame}
 	X_p = \frac{1}{|| \nabla H ||} K_1 \nabla H,\quad Y_p = \frac{1}{|| \nabla H ||} K_2 \nabla H
 \end{align}
satisfy $\omega(X_p,Y_p)=1$, and form a symplectic basis to the symplectic complement (\ref{sym_complement}), which is transverse to the Hamiltonian vector field $X_H = J \nabla H$ on an energy level set.\

Furthermore, the path (\ref{path}) decomposes symplectically into planar and spatial blocks:
\begin{align*}
	\Psi_{\gamma}(t) = \begin{pmatrix}
		\Psi_{\gamma}^p(t) & 0\\
		0 & \Psi_{\gamma}^s(t)
	\end{pmatrix},\quad \Psi_{\gamma}^p(t), \Psi_{\gamma}^s(t) \colon [0,T] \rightarrow Sp(2) = SL(2,\mathbb{R}).
\end{align*}
In view of the product property of the Maslov index (\ref{index_product}), the Conley--Zehnder index of the planar orbit~$\gamma$ splits additively into the planar and spatial index, which we denote by
$$ \mu_{CZ}(\gamma) = \mu_{CZ}^p (\gamma) + \mu_{CZ}^s (\gamma). $$
If $\gamma$ is non-degenerate, the definition of each index coincides with the construction by Hofer--Wysocki--Zehnder \cite{hofer_w_z}:\ Let $\gamma^k$ be the $k$-covering of $\gamma$ and assume that $\gamma^k$ is non-degenerate for all $k \geqslant 1$.\ The Conley--Zehnder index measures the rotation of corresponding eigenvectors, whose rotation angles are measured as a real number, and not modulo $2 \pi$, which is the case for rotation angles computed from Floquet multipliers.\ In the elliptic case, we denote by $\varphi_p$ and~$\varphi_s$ the corresponding rotation angles as a real number for the underlying single-turn orbit~$\gamma$.\ Then, each index is determined by
\begin{align*}
	\mu_{CZ}^p (\gamma^k) = 1 + 2 \cdot \lfloor k \cdot \varphi_p /(2\pi) \rfloor, \quad \mu_{CZ}^s (\gamma^k) = 1 + 2 \cdot \lfloor k \cdot \varphi_s /(2\pi) \rfloor.
\end{align*}
This means that each index measures the number of times that eigenvalues crosses 1 along the whole orbit $\gamma^k$, that jumps by~2 whenever 1 is crossed and is odd.\ It is worth to mention that in \cite{aydin_babylonian} we have used this Conley--Zehnder index construction to compute numerically the Babylonian lunar periods that are known from astronomical observations.\

In the hyperbolic case, the corresponding eigenvectors are rotated by $n\pi$ for an integer~$n$, and the corresponding index equals
$$ \mu_{CZ}^{p/s}(\gamma^k) = k n,\quad n \in \begin{cases}
	2 \mathbb{Z} & \text{ for the pos. hyperbolic case}\\
	2 \mathbb{Z} + 1 & \text{ for the neg. hyperbolic case.}
\end{cases} $$

\subsection{Interaction at bifurcation points}

\begin{figure}[t!]
	\centering
	\includegraphics[scale=1]{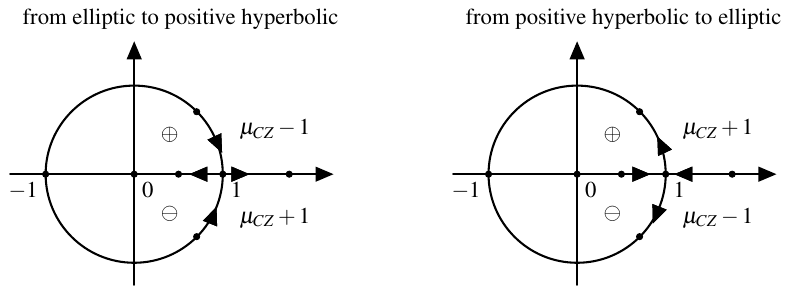}
	\caption{The index jump (reproduced from \cite{aydin_cz} with minor modifications).\ Left:\ When eigenvalue 1 is crossed from above (or below), the index goes down (or up) by 1.\ Right:\ When Floquet multipliers after crossing eigenvalue 1 goes up (or down), the index goes up (or down) by 1.\ Direction of crossing is determined by $\pm$ Krein signature.}
	\label{figure_3_index_jump}
\end{figure}

Along the continuation of an orbit family, when reaching a bifurcation point within a transition from elliptic to positive hyperbolic, or vice versa, the corresponding index jumps by $\pm 1$.\ This index jump depends on the direction of crossing the eigenvalue 1, as shown in Figure \ref{figure_3_index_jump}.\ When there is a touching of the eigenvalue~1, the index jumps by $\pm 2$, which is determined by the direction of the rotation angles.\ In order to determine the direction of the rotation, we consider the \textbf{Krein signature} \cite[Appendix~29]{arnold_avez} which relates a $\pm$ sign to each pair of elliptic Floquet multipliers:\ Let $M \in Sp(2n)$, and let $\lambda = e^{i \theta}$ and $\overline{\lambda} = e^{-i \theta}$ be a pair of elliptic eigenvalues of $M$ with $\lambda \overline{\lambda} \neq 1$.\ On the real 2-dimensional plane $P := \langle \text{Re}(v), \text{Im}(v) \rangle $, where $v$ is an eigenvector of~$M$ with eigenvalue $\lambda$, a quadratic form is defined by
\begin{align}\label{quadratic_form}
	Q(\xi,\xi) := \omega \left( \xi, M \xi \right),\quad \xi \in P,
\end{align}
which in matrix form is represented by
\begin{align}\label{quadratic_krein}
	Q = B^T J^T M B,
\end{align}
with $B = [ \text{Re}(v) \text{ } \text{Im}(v) ] \in \text{Mat}(2n \times 2,\mathbb{R})$.\ This quadratic form is non-degenerate, either positive or negative definite.\ The \textbf{Krein signature} associated to the elliptic pair $\{\lambda,\overline{\lambda}\}$ is defined as the signature of the quadratic form $Q$, which is a symplectic invariant due to its independent choice of the eigenvector.\ In particular, if the signature is positive then the rotation is determined by $\theta \in (0,\pi)$, and if the signature is negative then the rotation is determined by $-\theta$ in the range $(\pi,2\pi)$.\ Krein's main result was that the collision of two pairs of elliptic Floquet multipliers with different signatures on the unit circle can provoke complex instability, while the collision with same signatures does not cause a move off the unit circle.\ Fur our purpose, we use the Krein signature to specify the index jump, as illustrated in Figure~\ref{figure_3_index_jump}.\ Moreover, for VSR orbits the Krein signature determines the value of the integer $d$ in (\ref{vsr}).\ 

When working locally near a family of non-degenerate periodic orbits, then there is a fascinating topological bifurcation invariant, based on holomorphic curves \cite{ginzburg}:\ the \textbf{local Floer homology} and its \textbf{Euler characteristic}, the alternating sum of the ranks of the homology groups.\ The periodic orbits are the generators of the local homology groups that are graded by the Conley--Zehnder index.\ Therefore, the index provides important information how different families are related to each other at branch points.\ In practice, we check the invariance of the Euler characteristic by counting the periodic orbits with indices locally at bifurcation points.\ Furthermore, in the sense of SFT \cite{eliashberg_givental_hofer}, for the planar orbits $\gamma$, if 
\begin{align} \label{index_parity}
	\mu_{CZ}^p(\gamma^k) \equiv \mu_{CZ}^p(\gamma) \text{ mod } 2, \quad \mu_{CZ}^s(\gamma^k) \equiv \mu_{CZ}^s(\gamma) \text{ mod } 2,
\end{align}
or both equations in \eqref{index_parity} are not simultaneously satisfied, then $\gamma^k$ is referred to as \textit{good orbit}.\ Otherwise, $\gamma^k$ is called \textit{bad orbit} that are ignored in the local Floer homology.\ Therefore, only the good orbits will be counted with indices.\ Notice that all the single-turn periodic orbits are good, and the bad orbits occur as the $k$-covering of the orbits with exactly one pair of negative hyperbolic Floquet multipliers, where~$k$ is even.\ For the spatial solutions, we treat each pair of the Floquet multipliers and its interaction with the Conley--Zehnder index at bifurcation points in the same way as described before.\

\section{Analytical approach to comet-type periodic orbits}
\label{sec4}

In this section we provide an analytical approach to the two one-parameter families of planar nearly circular periodic orbits of the CR3BP which are located far away from the barycenter of the primaries.\ These solutions are the \textit{comet-type periodic orbits}.\ We first prove their existence by using the classical Poincaré continuation method, such as proposed in \cite[Chapter~9]{meyer_hall_offin}.\ Then, we determine their Conley--Zehnder indices by applying the construction from the previous section.\ The results in this section hold for all mass ratios of the CR3BP.\

\begin{theorem1}\label{theorem1}
	\textit{There exist two one-parameter families of planar nearly circular periodic orbits of the CR3BP which are at a great distance from the barycenter of the primaries.\ While one family exists for high energy values, the other family exists for low energy values.\ We denote these families by $\kappa_-$ and~$\kappa_+$, correspondingly.\ Moreover, the orbits of both families are of elliptic type and their Conley--Zehnder indices are given by}
	\begin{align}\label{comet_index}
		\mu_{CZ}(\kappa_{\pm}) = \mu_{CZ}^p(\kappa_{\pm}) + \mu_{CZ}^s(\kappa_{\pm}) = 2 = 1 + 1,
	\end{align}
	\textit{with positive planar and spatial Krein signatures.}
\end{theorem1}

The periodic solutions of the family $\kappa_+$ we refer to as the \textit{direct/prograde comet-type periodic orbits}, and those of the family $\kappa_-$ we refer to as the \textit{retrograde comet-type periodic orbits}.\ Before we prove Theorem \ref{theorem1}, we first discuss in the next theorem the Hamiltonian (\ref{hamiltonian}) of the CR3BP when we introduce a small scale parameter $\varepsilon$ to approximate such comet orbits with very large radii.\

\begin{theorem1}\label{theorem2}
	\textit{Let $\varepsilon > 0$ be a small parameter.\ After the conformally symplectic scaling}
	\begin{align}\label{scaling}
		(x,y,z,p_x,p_y,p_z) \mapsto \left( \varepsilon^{-2} x , \varepsilon^{-2} y, \varepsilon^{-2} z, \varepsilon p_x, \varepsilon p_y, \varepsilon p_z \right),
	\end{align}
	\textit{the Hamiltonian of the CR3BP (\ref{hamiltonian}) becomes}
	\begin{align}\label{hamiltonian_2}
		H = p_x y - p_y x + \varepsilon^3 \left( \frac{1}{2} \left( p_x^2 + p_y^2 + p_z^2 \right) - \frac{1}{(x^2 + y^2 + z^2)^{1/2}} \right) + \mathit{O}(\varepsilon^5).
	\end{align}
\end{theorem1}

From (\ref{hamiltonian_2}) we observe that when $\varepsilon$ is very small, which in view of the scaling (\ref{scaling}) means that the infinitesimal body is at a great distance from the primaries, the dominating force is the Coriolis force, and the next most significant force looks like a Kepler problem with both point masses at the origin.\

\begin{proof}[Proof of Theorem \ref{theorem2}]
	We first note that the scaling (\ref{scaling}) is conformally symplectic with the constant conformal factor $\varepsilon^{-1}$.\ To incorporate this magnification we consider a transformation of the Hamiltonian~(\ref{hamiltonian}) that is given by $H \rightarrow \varepsilon H \left( \varepsilon^{-2} x , \varepsilon^{-2} y, \varepsilon^{-2} z, \varepsilon p_x, \varepsilon p_y, \varepsilon p_z \right)$, which reads
	\begin{align}\label{hamiltonian_3}
		p_x y - p_y x + \frac{1}{2} \varepsilon^3 \left( p_x^2 + p_y^2 + p_z^2 \right) - \varepsilon \left( U_1^{\varepsilon}(x,y,z) + U_2^{\varepsilon}(x,y,z) \right),
	\end{align}
	where
	\begin{align}\label{func}
		U_1^{\varepsilon}(x,y,z) = \frac{1 - \mu}{\left(\tilde{r} + \mu^2 + 2\mu x/\varepsilon^2 \right)^{1/2}},\quad U_2^{\varepsilon}(x,y,z) = \frac{\mu}{\left(\tilde{r} + (1-\mu)^2 - 2(1-\mu)x/\varepsilon^2 \right)^{1/2}},
	\end{align}
	and $\tilde{r} := \left( x^2 + y^2 + z^2\right)/ \varepsilon^4 $.\ We next simplify (\ref{func}) by the Legendre expansions.\ To this end, we write
	\begin{align}\label{func_2}
		U_1^{\varepsilon}(x,y,z) = \varepsilon^2 \frac{1 - \mu}{r \left( 1 - 2 \bold x \bold t + \bold t^2 \right)^{1/2}},\quad U_2^{\varepsilon}(x,y,z) = \varepsilon^2 \frac{\mu}{r \left( 1 - 2 \tilde{ \bold x} \tilde{\bold t} + \tilde{ \bold t}^2 \right)^{1/2}},
	\end{align}
	where
	\begin{align}\label{func2_5}
		\bold x = - x/ r,\quad \bold t = \varepsilon^2 \mu / r,\quad \tilde{\bold x} = x/ r,\quad \tilde{\bold t} = \varepsilon^2 (1 - \mu)/ r, \quad r = \left( x^2 + y^2 + z^2 \right)^{1/2}.
	\end{align}
	Since $|\bold x| < 1$ and $|\tilde{\bold x}| < 1$, the corresponding expansions of (\ref{func_2}) in terms of the Legendre polynomials are given by
	\begin{align}\label{func_3}
		U_1^{\varepsilon}(x,y,z) = \varepsilon^2 \frac{1 - \mu}{r} \sum_{n = 0}^{\infty} P_n (\bold x) \bold t^n,\quad P_0 (\bold x) = 1,\quad P_1 (\bold x) = \bold x,
	\end{align}
	and
	\begin{align}\label{func_4}
		U_2^{\varepsilon}(x,y,z) = \varepsilon^2 \frac{\mu}{r} \sum_{m = 0}^{\infty} P_m (\tilde{\bold x}) \tilde {\bold t}^m,\quad P_0 (\tilde{\bold x}) = 1,\quad P_1 (\tilde{\bold x}) = \tilde{\bold x}.
	\end{align}
	By using the expressions (\ref{func_3}) and (\ref{func_4}), together with (\ref{func2_5}), we obtain
	\begin{align*}
		U_1^{\varepsilon}(x,y,z) + U_2^{\varepsilon}(x,y,z) = \varepsilon^2 \left( \frac{1}{r} - \frac{x(1 - 2 \mu)}{r^2} + \frac{1 - \mu}{r} \sum_{n = 2}^{\infty} P_n (\bold x) \bold t^n + \frac{\mu}{r} \sum_{m = 2}^{\infty} P_m (\tilde{\bold x}) \tilde {\bold t}^m \right),
	\end{align*}
	which, when inserted into (\ref{hamiltonian_3}), proves the theorem.\
\end{proof}

\begin{proof}[Proof of Theorem \ref{theorem1}]
	In cylindrical symplectic coordinates $(r,\theta,z,p_r,p_{\theta},p_z)$, introduced as follows,
	\begin{alignat*}{2}
		&x = r \cos \theta,\quad && \quad p_x = p_r \cos \theta - \frac{p_{\theta}}{r} \sin \theta,\\
		&y = r \sin \theta,\quad && \quad p_y = p_r \sin \theta + \frac{p_{\theta}}{r} \cos \theta,\\
		&z = z,\quad && \quad p_z = p_z,
	\end{alignat*}
	the Hamiltonian (\ref{hamiltonian_2}) becomes
	\begin{align}\label{ham_cyl}
		H = - p_{\theta} + \varepsilon^3 \left( \frac{1}{2} \left( p_r^2 + p_{\theta}^2 / r^2 + p_z^2 \right) - \frac{1}{(r^2 + z^2)^{1/2}} \right) + \mathit{O}(\varepsilon^5).
	\end{align}
	The equations of motion are given by
	\begin{alignat}{2}\label{equations}
		&\dot{r} = \varepsilon^3 p_r,\quad && \quad \dot{p}_r = \varepsilon^3 \left( \frac{p_{\theta}^2}{r^3} - \frac{r}{(r^2 + z^2)^{3/2}} \right),\nonumber\\
		&\dot{\theta} = \varepsilon^3 \frac{p_{\theta}}{r^2} - 1,\quad && \quad \dot{p}_{\theta} = 0,\\
		&\dot{z} = \varepsilon^3 p_z,\quad && \quad \dot{p}_z = - \varepsilon^3 \frac{z}{(r^2 + z^2)^{3/2}},\nonumber
	\end{alignat}
	where the terms of order $\varepsilon^5$ have been omitted.\ Within this order of approximation, $p_{\theta}$ is an integral.\ We shall now look for circular periodic solutions of (\ref{equations}) on the ecliptic $\{ z = p_z = 0 \}$.\ We find that the system~(\ref{equations}) accepts a pair of circular periodic solutions that are determined by
	\begin{align}\label{solutions}
		r = 1,\quad p_r = 0,\quad p_{\theta} = \pm 1,
	\end{align}
	with period $T = 2 \pi /(1 \mp \varepsilon^3)$.\ In order to linearize the equations (\ref{equations}) near the solutions~(\ref{solutions}), we expand $(r + \Delta r, \theta + \Delta \theta, z + \Delta z, p_r + \Delta p_r, p_{\theta} + \Delta p_{\theta}, p_z + \Delta p_z)$ near (\ref{solutions}), which leads to the autonomous linearized equations that split into 4-dimensional planar and 2-dimensional spatial components:\
	\begin{comment}
	\begin{align}\label{lin_eq}
		\begin{pmatrix}
			\Delta \dot{r}\\
			\Delta \dot{\theta}\\
			\Delta \dot{z}\\
			\Delta \dot{p}_r\\
			\Delta \dot{p}_{\theta}\\
			\Delta \dot{p}_z
		\end{pmatrix} = \begin{pmatrix}
		0 & 0 & 0 & \varepsilon^3 & 0 & 0\\
		\mp 2 \varepsilon^3 & 0 & 0 & 0 & \varepsilon^3 & 0\\
		0 & 0 & 0 & 0 & 0 & \varepsilon^3\\
		- \varepsilon^3 & 0 & 0 & 0 & \pm 2 \varepsilon^3 & 0\\
		0 & 0 & 0 & 0 & 0 & 0\\
		0 & 0 & - \varepsilon^3 & 0 & 0 & 0
		\end{pmatrix} \begin{pmatrix}
		\Delta r\\
		\Delta \theta\\
		\Delta z\\
		\Delta p_r\\
		\Delta p_{\theta}\\
		\Delta p_z
		\end{pmatrix}.
	\end{align}
	\end{comment}
	\begin{align}\label{lin_eq}
		\begin{pmatrix}
			\Delta \dot{r}\\
			\Delta \dot{\theta}\\
			\Delta \dot{p}_r\\
			\Delta \dot{p}_{\theta}
		\end{pmatrix} = \begin{pmatrix}
			0 & 0 & \varepsilon^3 & 0\\
			\mp 2 \varepsilon^3 & 0 & 0 & \varepsilon^3\\
			- \varepsilon^3 & 0 & 0 & \pm 2 \varepsilon^3\\
			0 & 0 & 0 & 0 \\
		\end{pmatrix} \begin{pmatrix}
			\Delta r\\
			\Delta \theta\\
			\Delta p_r\\
			\Delta p_{\theta}\\
		\end{pmatrix},\quad \begin{pmatrix}
		\Delta \dot{z}\\
		\Delta \dot{p}_z
		\end{pmatrix} = \begin{pmatrix}
		0 & \varepsilon^3\\
		- \varepsilon^3 & 0 
		\end{pmatrix} \begin{pmatrix}
		\Delta z\\
		\Delta p_z
		\end{pmatrix}.
	\end{align}
	The Jordan normal form of the corresponding matrix exponential is given by
	$$ \text{diag} \left\{ \begin{pmatrix}
		1 & 1\\
		0 & 1
	\end{pmatrix}, \begin{pmatrix}
	e^{i \varepsilon^3 t} & 0\\
	0 & e^{-i \varepsilon^3 t}
	\end{pmatrix} \right\},\quad \begin{pmatrix}
	e^{i \varepsilon^3 t} & 0\\
	0 & e^{-i \varepsilon^3 t}
	\end{pmatrix}. $$
	\begin{comment}
	$$ \begin{pmatrix}
		0 & 1 & 0 & 0 & 0 & 0\\
		0 & 0 & 0 & 0 & 0 & 0\\
		0 & 0 & i \varepsilon^3 & 0 & 0 & 0\\
		0 & 0 & 0 & - i \varepsilon^3 & 0 & 0\\
		0 & 0 & 0 & 0 & i \varepsilon^3 & 0 \\
		0 & 0 & 0 & 0 & 0 & - i \varepsilon^3
	\end{pmatrix}. $$
	\end{comment}
	Notice that the top-left $2 \times 2$ block appears since when the energy (or the period) varies, then the period (or the energy) varies as well.\ Moreover, at the first return time $T = 2 \pi/(1 \mp \varepsilon^3)$ we obtain the planar and spatial Floquet multipliers, both of which are of the form
	\begin{align}\label{floq}
		e^{\pm i \varepsilon^3 2 \pi /(1\mp \varepsilon^3)} = 1 \pm \varepsilon^3 2 \pi i + \mathit{O}(\varepsilon^6).
	\end{align}
	With respect to the planar Floquet multipliers, by using the Poincaré continuation method \cite[Chapter~9]{meyer_hall_offin}, we conclude that the two solutions (\ref{solutions}) of the equations (\ref{equations}) can be continued to the full equations related to (\ref{ham_cyl}), where the $\mathit{O}(\varepsilon^5)$ terms are included.\ This proves the assertion on the existence in the theorem.\
	
	We denote by $\kappa_+$ the family whose orbits are associated to $p_{\theta} = 1$, and by~$\kappa_-$ the family whose members are related to $p_{\theta} = -1$.\ For the computation of their Conley--Zehnder indices, we notice that the Hamiltonian vector field at the solutions (\ref{solutions}) is of the form $X_H = (\pm \varepsilon^3 - 1)\partial_{\theta}$.\ In view of the trivializations (\ref{spatial_frame}) and (\ref{planar_frame}), the planar and spatial frames $[X_p,Y_p]$ and $[X_s,Y_s]$ are formed by the vectors
	$$ X_p = \partial_{p_r},\quad Y_p = \partial_r,\quad X_s = \partial_{p_z},\quad Y_s = \partial_z. $$
	Taking into account (\ref{lin_eq}), the transverse linearized flow with respect to these frames is given by
	$$ \text{diag} \left\{ \begin{pmatrix}
		0 & -\varepsilon^3\\
		\varepsilon^3 & 0
	\end{pmatrix}, \begin{pmatrix}
	0 & -\varepsilon^3\\
	\varepsilon^3 & 0
	\end{pmatrix} \right\}, $$
	whose matrix exponential provides a path of symplectic matrices that is represented by
	\begin{align}\label{path_2}
		\Psi_{\kappa_{\pm}}(t) = \begin{pmatrix}
			\Psi_{\kappa_{\pm}}^p(t) & 0\\
			0 & \Psi_{\kappa_{\pm}}^s(t)
		\end{pmatrix} = \begin{pmatrix}
			\cos(t \varepsilon^3) & - \sin (t \varepsilon^3) & 0 & 0\\
			\sin (t \varepsilon^3) & \cos (t \varepsilon^3) & 0 & 0\\
			0 & 0 & \cos (t \varepsilon^3) & - \sin (t \varepsilon^3)\\
			0 & 0 & \sin (t \varepsilon^3)& \cos (t \varepsilon^3)
		\end{pmatrix},
	\end{align}
	which has crossings at $t \in 2 \pi \mathbb{Z} / \varepsilon^3$.\ Furthermore, the crossing form,
	$$ Q_t = \begin{pmatrix}
		0 & 1 & 0 & 0\\
		-1 & 0 & 0 & 0\\
		0 & 0 & 0 & 1\\
		0 & 0 & -1 & 0
	\end{pmatrix} \dot{\Psi}_{\kappa_{\pm}}(t) = \begin{pmatrix}
		\varepsilon^3 & 0 & 0 & 0\\
		0 & \varepsilon^3 & 0 & 0\\
		0 & 0 & \varepsilon^3 & 0\\
		0 & 0 & 0 & \varepsilon^3
	\end{pmatrix}, $$
	has always planar and spatial signature 2.\ We observe that the path (\ref{path_2}) rotates extremely slowly with a very large period of $2 \pi / \varepsilon^3$.\ It is evident that during one period of the solutions (\ref{solutions}) there is no crossing, which, in view of the definition of the Maslov index (\ref{maslov}), proves the Conley--Zehnder indices as stated in~(\ref{comet_index}).\ Furthermore, eigenvectors to the Floquet multipliers (\ref{floq}) are of the form $[\pm i \text{ } 1]^T$.\ By using (\ref{quadratic_krein}), it is now a straightforward computation that on the plane $P = \langle [0 \text{ } 1]^T, [1 \text{ } 0]^T \rangle$ the quadratic form (\ref{quadratic_form}) is positive definite, which means that the $\kappa_{\pm}$ orbits have planar and spatial positive Krein signatures.\
\end{proof}

\section{Numerical computational results within the Earth--Moon CR3BP}
\label{sec5}

For the results in this section, we have performed numerical integration of the equations of motion~(\ref{ham_equation}) and its variations~(\ref{var_eq}) by using an explicit Runge-Kutta method of order twelve, such as proposed in \cite{feagin}, which we have encoded in Python.\ All periodic solutions have been obtained via a classical corrector--predictor continuation method \cite{robin_markellos} (that is described briefly below).\ The mass ratio assumed for the Earth--Moon environment is $\mu = 1/82.27 \approx 0.0121550991$.\ We also assume that the unit of distance equals 384400 km, and the unit of time equals 4.348 days, that is defined as the sidereal lunar month (that is the orbital period of the Moon around the Earth $\approx$~27.32~days) divided by $2 \pi$.\ To compute the closest approach to the Earth and the lunar surface we further assume that the radius of the Earth and the Moon equals correspondingly 6371 km and 1737.4~km.\ We also analyze \textit{M:N resonant orbits}.\ These are the orbits that are in resonance with the Earth and the Moon, where $M$ corresponds to the number of revolutions on the orbit and $N$ to the number of revolutions of the Earth and the Moon.\ In particular, the~1:$k$ resonant orbits are those with period $k \cdot 2 \pi$ that corresponds to $k$ times the sidereal lunar month.\

\subsection{Retrograde and direct comet-type periodic orbits}

A first decent initial guess is provided by large circular Keplerian motions centered at the origin.\ We therefore consider the initial data of planar circular periodic solutions of the rotating Kepler problem (which is a special case of the CR3BP when $\mu = 0$), that are of the form
\begin{align}\label{initial_guess}
	(x_0,0,0,0,\dot{y}_0,0),\quad \dot{y}_0 = \pm \frac{1}{\sqrt{x_0}} - x_0,
\end{align}
where $x_0 > 0$ is the Keplerian semi-major axis and $\pm$ indicates the respective direct and retrograde motion.\ The corresponding energy value equals $- 1/(2 x_0) \mp \sqrt{x_0}$.\ For large semi-major axis $x_0$, these orbits have the following properties:\
\begin{itemize}[noitemsep]
	\item the Keplerian period equals $T_{\pm} = 2 \pi /(1 \mp x_0^{-3/2})$ which is slightly above $2 \pi$ for the direct, and slightly below $2 \pi$ for the retrograde orbits;
	\item the energy values are very low for the direct and very high for the retrograde orbits, or equivalently, the Jacobi constant values are very high for the direct and very low for the retrograde orbits;
	\item while the retrograde orbits are retrograde in both the rotating and inertial frame, the direct orbits are direct in the inertial and retrograde in the rotating frame; and
	\item the orbits are simple symmetric with only two perpendicular intersections with the $x$-axis per period.
\end{itemize}

By choosing an initial guess of the form (\ref{initial_guess}) and integrating the equations of motion (\ref{ham_equation}), the symmetric orbit will be periodic if at $t = T/2$ it fulfills the periodicity conditions:
\begin{align}\label{periodicity_conditions}
	y_{cut} \equiv y(x_0,0,0,0,\dot{y}_0,0;T/2) = 0,\quad \dot{x}_{cut} \equiv \dot{x}(x_0,0,0,0,\dot{y}_0,0;T/2) = 0.
\end{align}
Since these conditions are only approximately satisfied, we seek appropriate corrections $\delta x_0, \delta \dot{y}_0, \delta t$ to fulfill these conditions with a specified accuracy.\ For this purpose we consider the first-order Taylor expansion of (\ref{periodicity_conditions}) and obtain
\begin{align}\label{first_order_Taylor}
	y_{cut} + v_{21} \delta x_0 + v_{25} \delta \dot{y}_0 + \dot{y}_{cut}(t) \delta t = 0,\\
	\dot{x}_{cut} + v_{41} \delta x_0 + v_{45} \delta \dot{y}_0 + \ddot{x}_{cut}(t) \delta t = 0,\nonumber
\end{align}
where $v_{21} = \partial y_{cut} / \partial x_0$, $v_{25} = \partial y_{cut} / \partial \dot{y}_0$, $v_{41} = \partial \dot{x}_{cut} / \partial x_0$ and $v_{45} = \partial \dot{x}_{cut} / \partial \dot{y}_0$ are the corresponding entries of the variational matrix $V(t) = (v_{ij}) \in \mathbb{R}^{6 \times 6}$ determined by (\ref{var_eq}).\ By fixing one of the three unknowns in (\ref{first_order_Taylor}), we solve (\ref{first_order_Taylor}) for the remaining two, and start a new integration of the equations (\ref{ham_equation}) with the corrected initial data.\ We repeat the procedure in an iterative way of correction steps until a specified accuracy is satisfied, that is measured in the norm $\sqrt{y_{cut}^2 + \dot{x}_{cut}^2}$ which in this work is less than~$10^{-8}$.\ When a symmetric periodic orbit has been determined, we consider an arbitrary step of the initial parameter that we have fixed in the corrector steps, and obtain for the continuation the predictor scheme as previously.\ The choice of which of the parameters to fix has an important effect on the convergence of the solutions.\ If the fixed parameter reaches an extremum, the algorithm will break off and a new choice has to be considered.\

\begin{figure}[t!]
	\centering
	\includegraphics[width=1\linewidth]{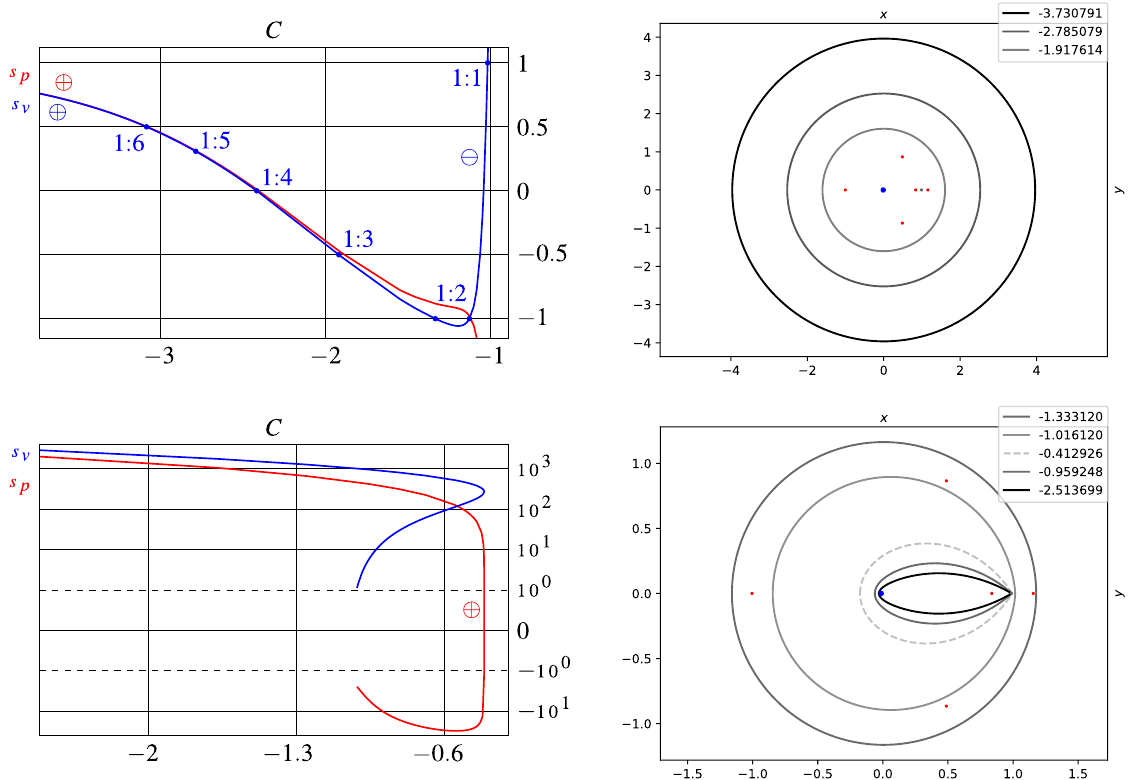}
	\caption{The family $\kappa_-$ of the retrograde comet-type periodic orbits.\ Right top shows the orbits that start at low Jacobi constants; from dark to light indicates increasing of Jacobi constants.\ Right bottom shows the continuation of the orbits; grey dashed is an orbit of birth-death type.\ Left top shows planar and vertical stability diagrams ($s_p$ and $s_v$) that are continued left bottom with a logarithmic scale.\ In the elliptic range at left top, the $d$:$k$ VSR orbits are denoted with respect to $s_v = \cos (2\pi d/k)$.\ The ``$\pm$'' sign indicates the Krein signature.}
	\label{plot_1}
\end{figure}

\begin{figure}[t!]
	\centering
	\includegraphics[width=1\linewidth]{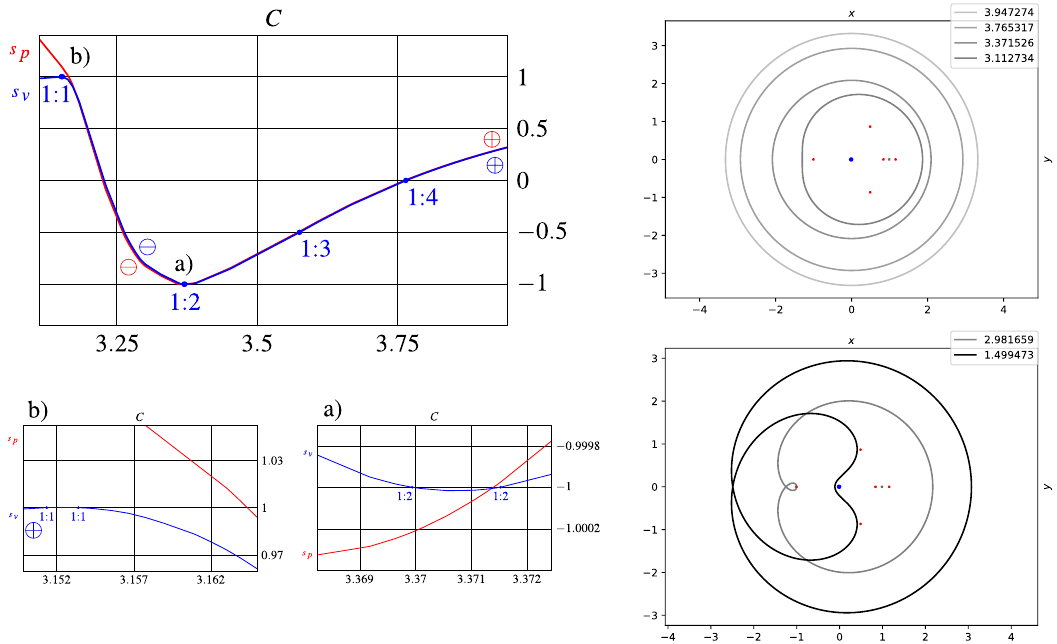}
	\caption{Family $\kappa_+$ of direct comet-type periodic orbits, continued up to approaching collision with the Earth.\ Right top shows $\kappa_+$ orbits starting at high Jacobi constants; from light to dark indicates decreasing of Jacobi constants.\ Right bottom shows continuation of $\kappa_+$ orbits.\ Left top shows planar and vertical stability diagrams ($s_p$ and $s_v$).\ Left bottom shows zoomed regions, a) with respect to crossings of $-1$, and b) with respect to crossings of 1.\ In the elliptic range at left top, the $d$:$k$ VSR orbits are denoted with respect to $s_v = \cos (2\pi d/k)$.\ The ``$\pm$'' sign indicates Krein signature.}
	\label{plot_2}
\end{figure}

Using the above corrector--predictor algorithm we have computed the two families $\kappa_{\pm}$ of the direct and retrograde comet-type periodic orbits.\ Figure \ref{plot_1} and Figure \ref{plot_2} shows orbit plots for certain Jacobi constant values, as well as the planar and vertical stability diagrams.\ In the elliptic region, the vertical self-resonant (VSR) orbits are identified that we denote in the form ``$d$:$k$'' with respect to the vertical stability index $s_v = \cos (2\pi d/k)$.\ The databases are collected in Table \ref{data_1}.\

\textit{Family $\kappa_-$ of retrograde comet orbits.}\ Along the continuation, the shape of the orbits smoothly changes from the circular type to an ellipse-type orbit around both primaries (see Figure \ref{plot_1}).\ The Jacobi constant values increase and reach a maximal value at $C \approx -0.412925$ that corresponds to a birth-death bifurcation (or saddle-node bifurcation which is a branch point from which two families bifurcate into the same energy direction with an index difference of 1.)\ Then, the Jacobi constant values decrease further to very low values.\ As it was already observed in \cite{broucke}, at the limit, there is a rectilinear orbit between both primaries with infinite velocity and infinite energy.\ Recall from Theorem \ref{theorem1} that the comet orbits in the generating procedure have Conley--Zehnder index $\mu_{CZ}(\kappa_-) = \mu_{CZ}^p(\kappa_-) + \mu_{CZ}^s(\kappa_-) = 2 = 1+1$, with positive Krein signature.\ From the stability diagrams in Figure~\ref{plot_1} we observe that both of the planar and vertical stability indices, $s_p$ and $s_v$, decrease, whereby $s_v$ is slightly below~$s_p$.\ At the Jacobi constant value $C \approx -1.333119$ the vertical stability index $s_v$ crosses $-1$, and after a short negative hyperbolic phase it becomes elliptic again at $C \approx -1.126762$, with negative Krein signature.\ After that, the vertical index increases extremely rapidly, crosses $+1$ at $C \approx -1.016119$ and grows strongly in the positive hyperbolic direction.\ At the latter transition, due to the negative Krein signature, the spatial Conley--Zehnder index jumps from 1 to 2.\ On the vertical stability diagram in Figure \ref{plot_1}, we have identified the first occurring VSR orbits up to the multiplicity six, that are labeled by 1:6, 1:5, 1:4, 1:3, 1:2 (two times) and 1:1.\ Shortly after the vertical stability index $s_v$ crosses $-1$ the second time, the planar stability index $s_p$ also crosses~$-1$.\ After that, $s_p$ remains negative hyperbolic for a while, increases very quickly, becomes elliptic with positive Krein signature, crosses $+1$ (birth-death bifurcation), and increases further in the positive hyperbolic way, while moving closely below the vertical stability index.\ At the birth-death branch point, since the Krein signature is positive, the planar Conley--Zehnder index jumps from 1 to 0.\

\textit{Family $\kappa_+$ of direct comet orbits.}\ Along the continuation, a symmetric point on the trajectory moves closer to the libration point $L_3$ and forces a loop, such that the trajectories only intersect the $x$-axis perpendicularly when they cross it for the second time (see Figure \ref{plot_2}).\ This loop then approaches collision with the Earth.\ We stop the continuation at $C \approx 1.499473$ where the orbit has the closest approach of 30974.12 km to the Earth surface.\ Recall from Theorem~\ref{theorem1} that the comet orbits in the generating procedure have Conley--Zehnder index $\mu_{CZ}(\kappa_+) = \mu_{CZ}^p(\kappa_+) + \mu_{CZ}^s(\kappa_+) = 2 = 1+1$, with positive Krein signature.\ As depicted in Figure \ref{plot_2}, the planar and vertical stability indices behave similarly, with the planar index being slightly above the vertical index.\ Both decrease, become negative hyperbolic for a very short phase, and then become elliptic again with negative Krein signature.\ After that, both stability indices cross $+1$.\ At this transition, since the Krein signatures are negative, both Conley--Zehnder indices jump from 1 to 2.\ Then, while the planar stability index remains positive hyperbolic, the vertical stability index becomes very quickly elliptic again with positive Krein signature, where the spatial Conley--Zehnder index jumps from 2 to 3.\ Furthermore, on the vertical stability diagram in Figure \ref{plot_2}, we have identified VSR orbits up to the multiplicity four, that are denoted by 1:4, 1:3, 1:2 (two times) and 1:1 (two times).

\subsection{Out-of-plane bifurcation results}

By adapting the corrector--predictor method (\ref{first_order_Taylor}) for spatial symmetric solutions, together with the computational technique of the Conley--Zehnder index discussed in Section \ref{sec3}, we have analyzed all the vertical bifurcations that emanate from the VSR orbits identified above.\ We illustrate our main results in the form of bifurcation graphs that provide a topological perspective, that are constructed in the following way.\ The vertices are branch points and the edges are the orbit families, labeled with Conley--Zehnder index.\ We draw from bottom to top in the direction of decreasing the Jacobi constant.\ Bold Jacobi constants correspond to critical branch points.\ Thick and think families indicate doubly and simple symmetric solutions, correspondingly.\ Edges at branch points associated to planar orbits, we draw in black, vertically and shortly before and after the bifurcations.\ Edges of vertical bifurcations we draw in red or blue:\ if the bifurcated branch begins with the same index as the underlying comet orbits before the bifurcation, then the bifurcated branch is referred to as  the\textit{red branch}, otherwise it is referred to as the \textit{blue branch}.\ Symmetric branches that are obtained by the $\sigma$-symmetry (reflection at the ecliptic) have same color and are dashed.\ Families labeled with overlined index indicate bad orbits.\ By ``b-d'' we denote birth-death bifurcation, whose local Floer homology and its Euler characteristic are zero.\

\subsubsection{Single-turn and period-doubling bifurcations from the retrograde comet orbits}

\begin{figure}[t!]
	\centering
	\includegraphics[width=1\linewidth]{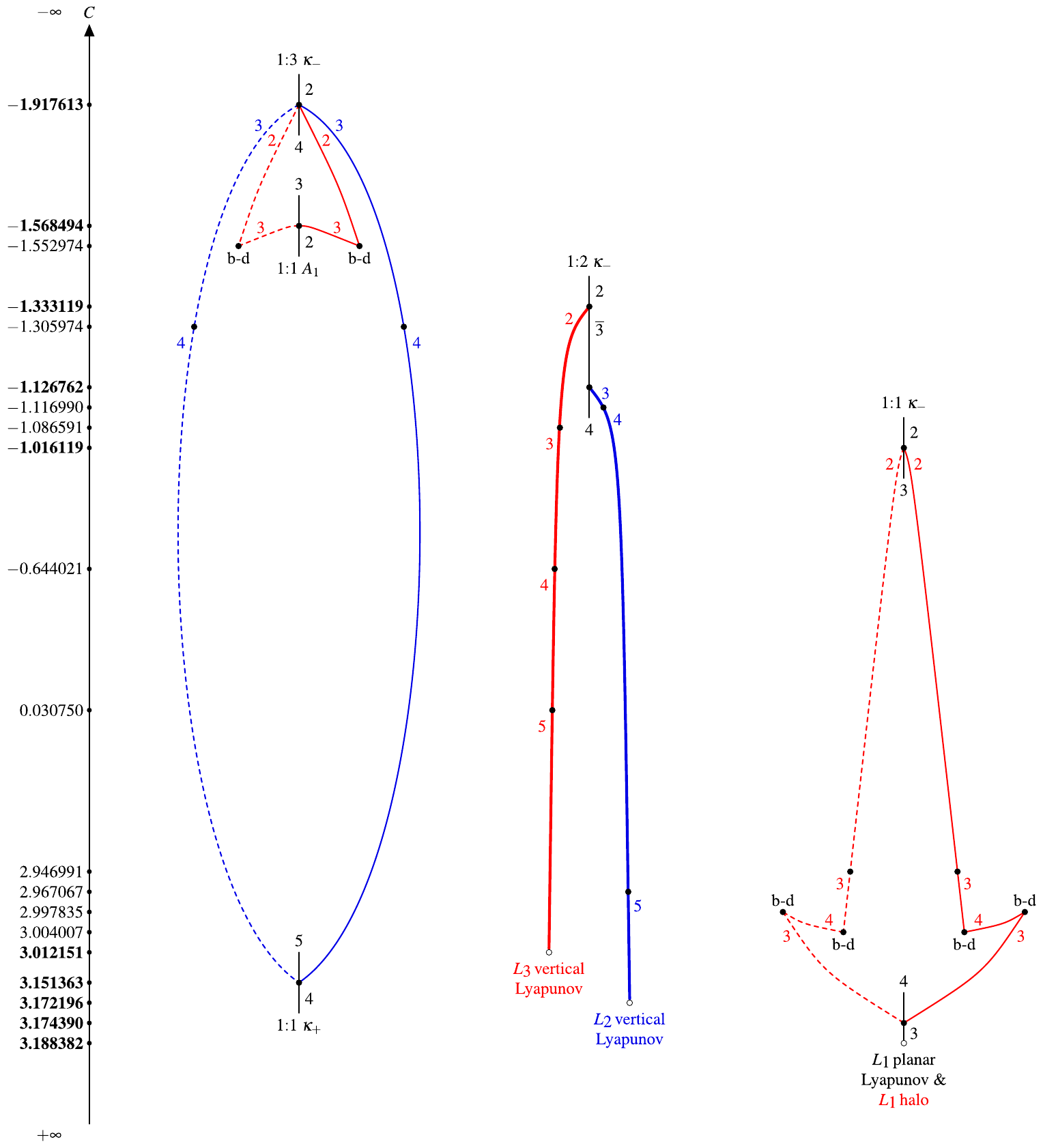}
	\caption{Bifurcation diagram associated to bifurcated branches from 1:1, 1:2 and 1:3 retrograde $\kappa_-$ orbits.\ The bifurcated branches from 1:1 and 1:2 retrograde $\kappa_-$ orbits correspond to $L_1$ halo, $L_2$ and $L_3$ vertical Lyapunov orbit families.\ The bifurcated branches from 1:3 retrograde $\kappa_-$ orbits form bridge families to 1:1 $A1$ and to 1:1 direct $\kappa_+$ orbits.}
	\label{plot_bif_1}
\end{figure}

\begin{figure}[t!]
	\centering
	\includegraphics[width=1\linewidth]{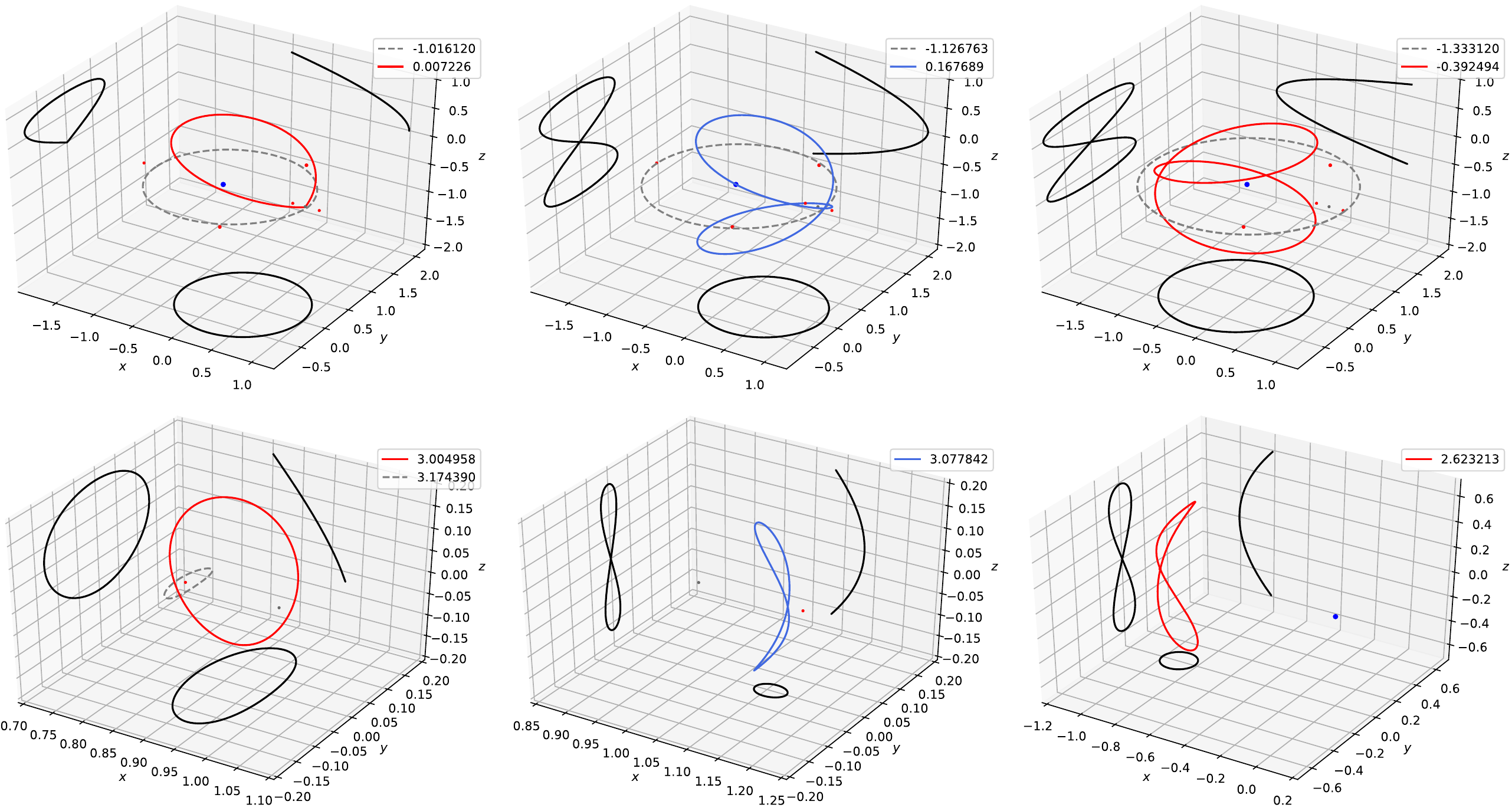}
	\caption{Left, middle and right show orbit plots of bifurcated branches from 1:1, 1:2 and 1:2 retrograde comet $\kappa_-$ orbits (grey dashed top), that correspond to the families of $L_1$ halo north (terminating at $L_1$ planar Lyapunov; grey dashed bottom left), $L_2$ and $L_3$ vertical Lyapunov orbits (figure-eight shape), respectively.\ Red orbits on the left are symmetric w.r.t.\ the $xz$-plane, while blue and red orbits in the middle and on the right are doubly symmetric w.r.t.\ the $xz$-plane and the $x$-axis.}
	\label{plot_3}
\end{figure}

The out-of-plane bifurcated branches from the 1:1 and 1:2 retrograde $\kappa_-$ orbits are well-known, but usually generated from the collinear libration points $L_i$ $(i=1,2,3)$:\ The one that bifurcates from the 1:1~retrograde $\kappa_-$ orbits corresponds to the $L_1$ halo family, and the ones that bifurcate from the 1:2~retrograde $\kappa_-$ orbits correspond to the $L_2$ and $L_3$ vertical Lyapunov families.\ The bifurcation diagram that illustrates these bifurcation results is provided in Figure \ref{plot_bif_1}, some orbit plots are shown in Figure \ref{plot_3}, and database are collected in Table \ref{data_2}.\ While the $L_2$ and $L_3$ vertical Lyapunov orbits are doubly symmetric w.r.t.\ the $xz$-plane and the $x$-axis, the $L_1$ halo orbits are simple symmetric w.r.t.\ the $xz$-plane.\ By using the $\sigma$ symmetry, the $L_1$ halo family is commonly categorized into \textit{northern} and \textit{southern halo family}, according to the orientation of the orbits relative to the ecliptic.\ While the southern halo orbits extend below the ecliptic, the northern halo orbits extend above the ecliptic.\ Since these orbit families are well-known, we only discuss that the Conley--Zehnder indices associated to these bifurcations are in accordance with the Euler characteristic.\ In view of Figure~\ref{plot_bif_1}, in the case of the 1:1 retrograde $\kappa_-$ orbits, before the branch point there are only the retrograde $\kappa_-$ orbits with index 2, and after the branch point the retrograde $\kappa_-$ orbits have index 3 and each of the two $L_1$ halo branches has index 2.\ Therefore, the Euler characteristic before and after the bifurcation point is correspondingly given by
$$(-1)^2 = 1,\quad 2\cdot (-1)^2 + (-1)^3 = 1.$$
Between the two period-doubling bifurcation points of the 1:2 retrograde $\kappa_-$ orbits, the orbits are vertical negative hyperbolic and planar elliptic.\ Therefore, all the double-covered retrograde $\kappa_-$ orbits between these two period-doubling bifurcation points are bad orbits with index $\overline{3}$, as labeled in Figure \ref{plot_bif_1}.\ Before and after the first branch point of the 1:2 retrograde $\kappa_-$ orbits, the double-covered orbits have index~2, while the bifurcated $L_3$ vertical Lyapunov orbits have index~2.\ Therefore, the Euler characteristic before and after the bifurcation is $(-1)^2 = 1$.\ Before the bifurcation at the second branch point of the 1:2 retrograde $\kappa_-$ orbits, there are only the double-covered orbits that are bad with index $\overline{3}$, while after the bifurcation there are the double-covered orbits with index 4 and the bifurcated $L_2$ vertical Lyapunov orbits with index 3.\ Therefore, the Euler characteristic before and after the bifurcation point equals zero.\

\begin{figure}[t!]
	\centering
	\includegraphics[width=1\linewidth]{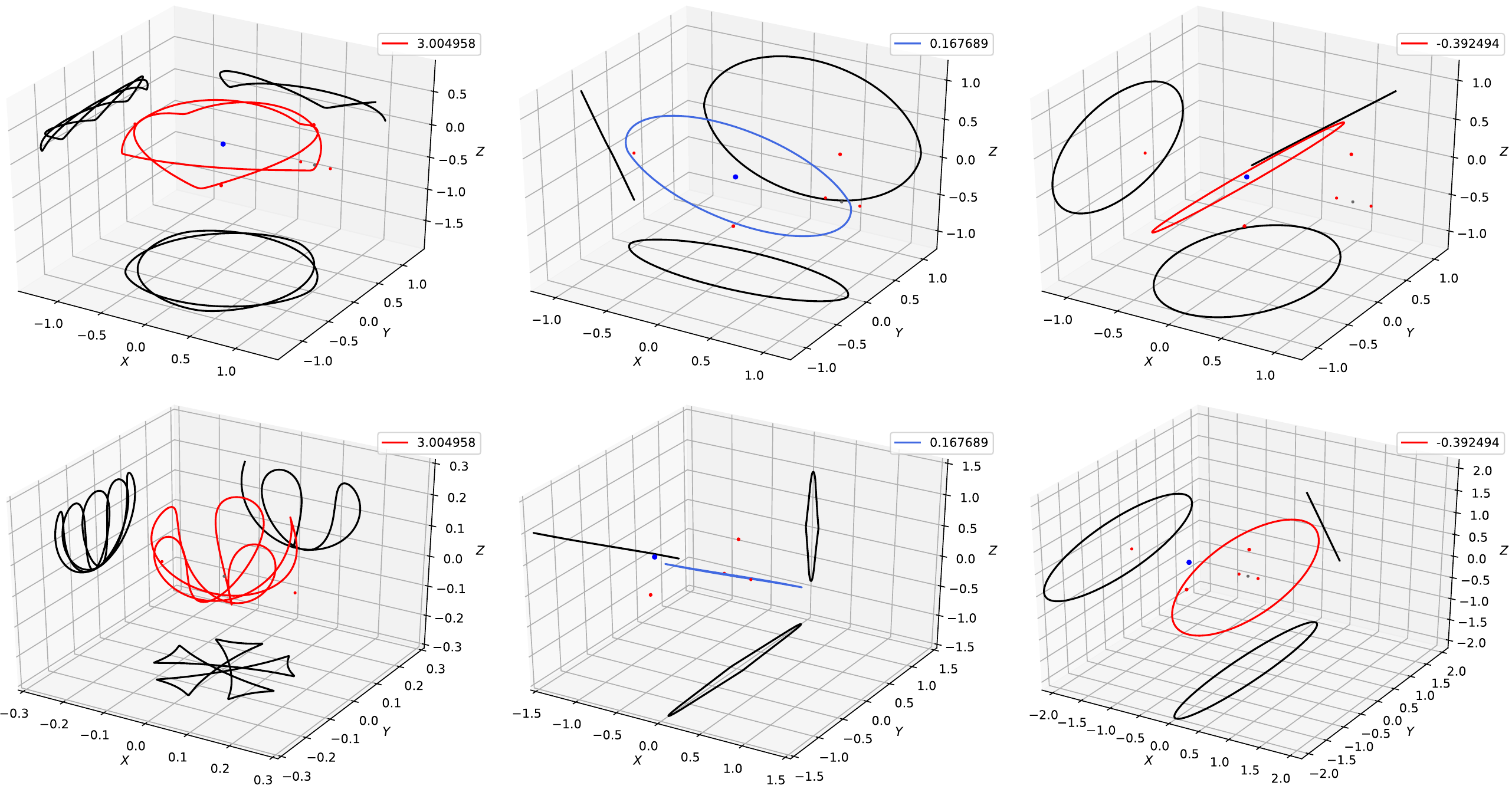}
	\caption{Left, middle and right show 5:2, 1:1 and 1:1 resonant orbits of $L_1$ halo, $L_2$ and $L_3$ vertical Lyapunov orbit families, respectively.\ Plots at the top are in the inertial frame, and plots at the bottom are in the Moon-centered inertial frame, after five and one revolutions of the Earth and the Moon, respectively.\ Corresponding orbits are plotted in the rotating frame in Figure \ref{plot_3}.}
	\label{plot_3_2}
\end{figure}

Along the continuation of the bifurcated branches, we also find resonant orbits.\ Within the $L_1$ halo orbits, the orbit at $C \approx 3.004958$ is a 5:2 resonant orbit with period of 10.92 days, and the closest approaches of 25974.84 km and 333024.03 km to the lunar and Earth surface, correspondingly.\ The orbit is plotted in Figure \ref{plot_3_2}, which has elliptic and positive hyperbolic Floquet multipliers of the form $\{ e^{\pm i \theta}, \lambda, \lambda^{-1} \}$, with $\theta \approx 2.438$ and $\lambda \approx 17.138$.\ Along each of the $L_2$ and $L_3$ vertical Lyapunov orbits there are 1:1 resonant orbits that are plotted in Figure \ref{plot_3_2}.\ Both are positive hyperbolic (see Table \ref{data_2}).\

\subsubsection{Bifurcations from the $k$-covered retrograde comet orbits, for $k=3,4,5,6$}

\begin{figure}[t!]
	\centering
	\includegraphics[width=1\linewidth]{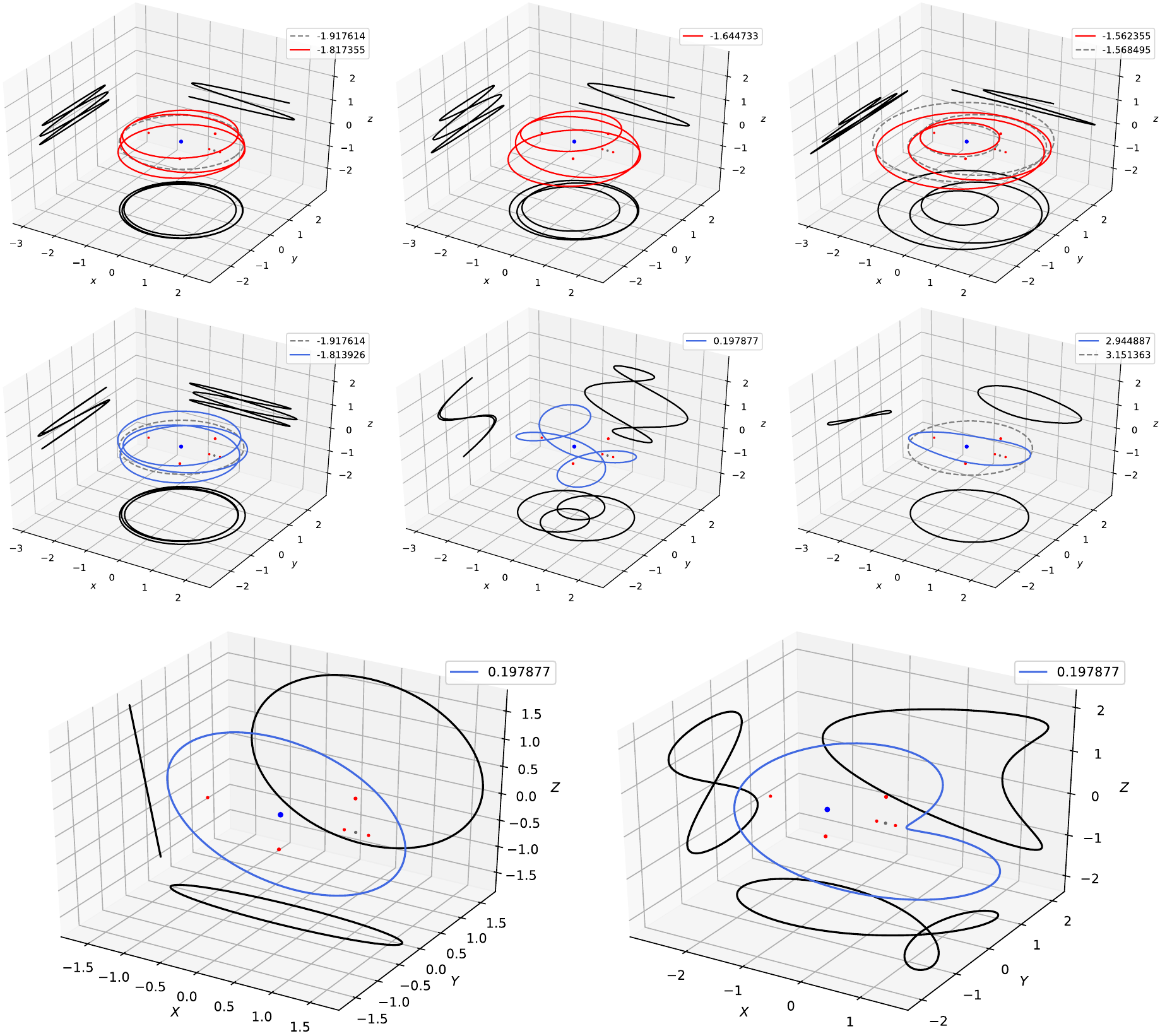}
	\caption{Top and middle:\ Red and blue orbits that bifurcate from 1:3 $\kappa_-$ (grey dashed on the left), continued from left to right.\ While the red orbits are symmetric w.r.t.\ the $xz$-plane and terminate at 1:1~$A_1$ (grey dashed on the right), the blue orbits are symmetric w.r.t.\ the $x$-axis and terminate at 1:1 $\kappa_+$ (grey dashed on the right).\ Middle shows the 1:2 resonant blue orbit in the rotating frame, that is plotted at the bottom in the inertial frame (left) and in the Moon-centered inertial frame (right), after two revolutions of the Earth and the Moon.}
	\label{plot_4}
\end{figure}

\setlength{\fboxsep}{1pt}
\framebox{\textit{(1) Bifurcation from 1:3 retrograde $\kappa_-$ orbit.}}\ At the Jacobi constant value $C \approx -1.917613$ the Conley--Zehnder index of the 3-covering of the 1:3 retrograde $\kappa_-$ orbit with period $T \approx 12.619274$ ($\approx$ 54.86~days) jumps from 2 to~4.\ This orbit has the closest approaches of 235751.76 km and 605166.71 km to the lunar and Earth surface, respectively.\ The members of the red branch start with index 2 and are simple symmetric w.r.t.\ the $xz$-plane, while the orbits of the blue branch start with index 3 and are simple symmetric w.r.t.\ the $x$-axis.\ Each of the two branches has a symmetric branch by using the $\sigma$-symmetry, and their continuation terminates at the ecliptic, as illustrated in the bifurcation diagram in Figure \ref{plot_bif_1}.\ Some orbits are plotted in Figure~\ref{plot_4}, and data are provided in Table \ref{data_3}.\ The out-of-plane initial parameter of the red orbits reaches the maximal value $z \approx 0.535093$ at $C \approx -1.644732$, and then, after a birth-death bifurcation at $C \approx -1.552974$, the red branch terminates planar at $C \approx -1.568494$ on a vertical branch point that belongs to the family of retrograde periodic solutions around the Earth, denoted by $A_1$ in Broucke's notation \cite{broucke}.\ All periodic orbits of the red family with Conley--Zehnder index 2 are stable.\ The blue branch grows until to $\dot{z} \approx 0.808557$ at $C \approx 0.845579$.\ Then, it goes down to plane on the second critical 1:1 direct $\kappa_+$ orbit at $C \approx 3.151363$ with period $T \approx 12.547041$ ($\approx$ 54.55 days) and the minimal distances of 271153.26 km and 567450.19 km to the lunar and Earth surface, respectively.\ This is our first result of a bridge family between the retrograde and direct comet orbits that consists of spatial periodic solutions.\ Along the continuation of the blue branch, some loops appear and disappear on the orbits, as can be seen in Figure \ref{plot_4}.\ All the blue orbits are unstable.\ Furthermore, within the blue family at $C \approx 0.197877$ we find a 1:2 resonant orbit (shown in Figure \ref{plot_4}) with index 4 and positive hyperbolic Floquet multipliers close to 1 (see Table \ref{data_3}).\

\begin{figure}[t!]
	\centering
	\includegraphics[width=1\linewidth]{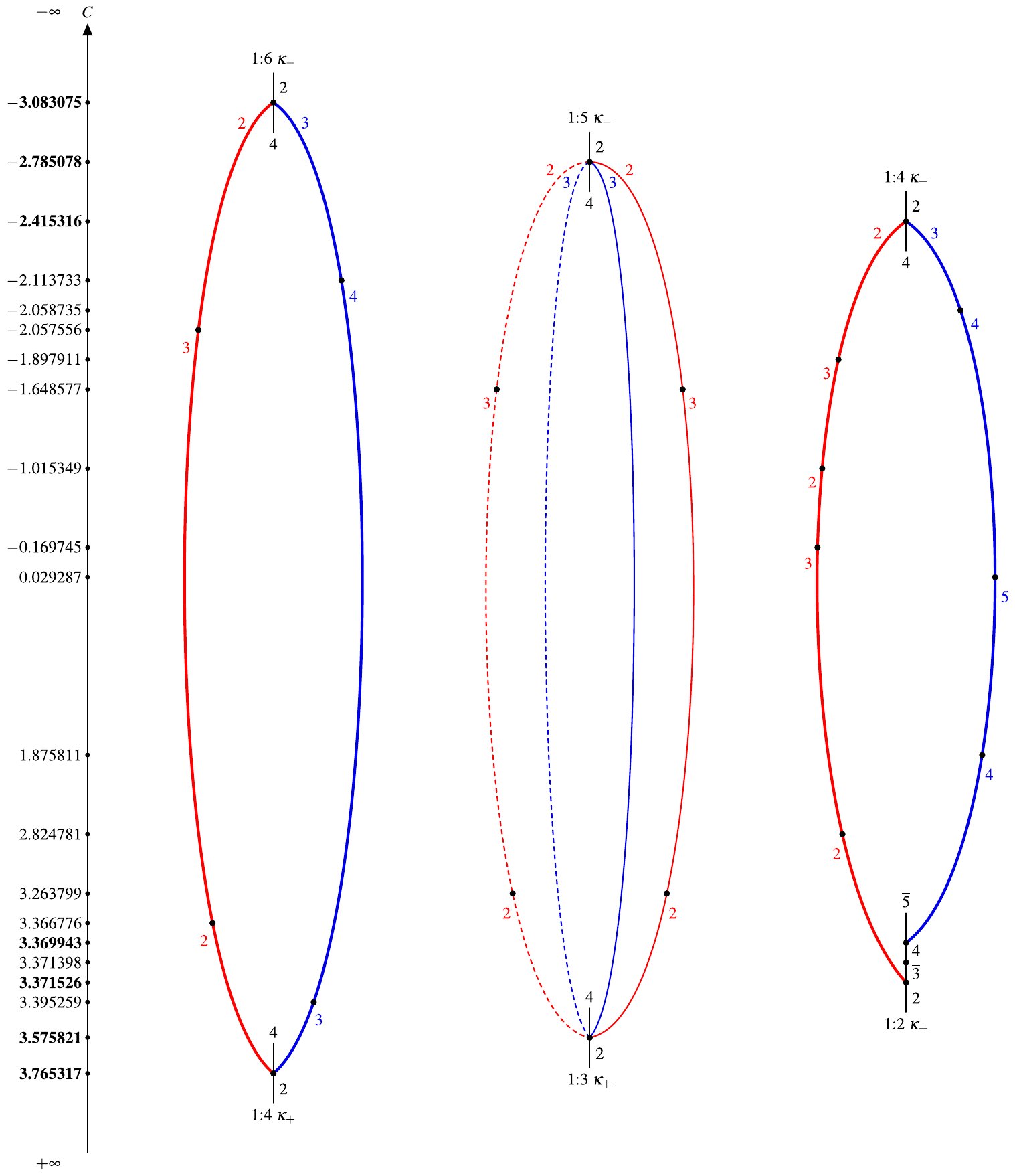}
	\caption{Bifurcation graph related to the bridge families between the 1:$k$ retrograde comet $\kappa_-$ orbits and the 1:$\tilde{k}$ direct comet $\kappa_+$ orbits, with $k=4,5,6$, and $\tilde{k} = k-2$.}
	\label{plot_bif_2}
\end{figure}

\begin{figure}[t!]
	\centering
	\includegraphics[width=1\linewidth]{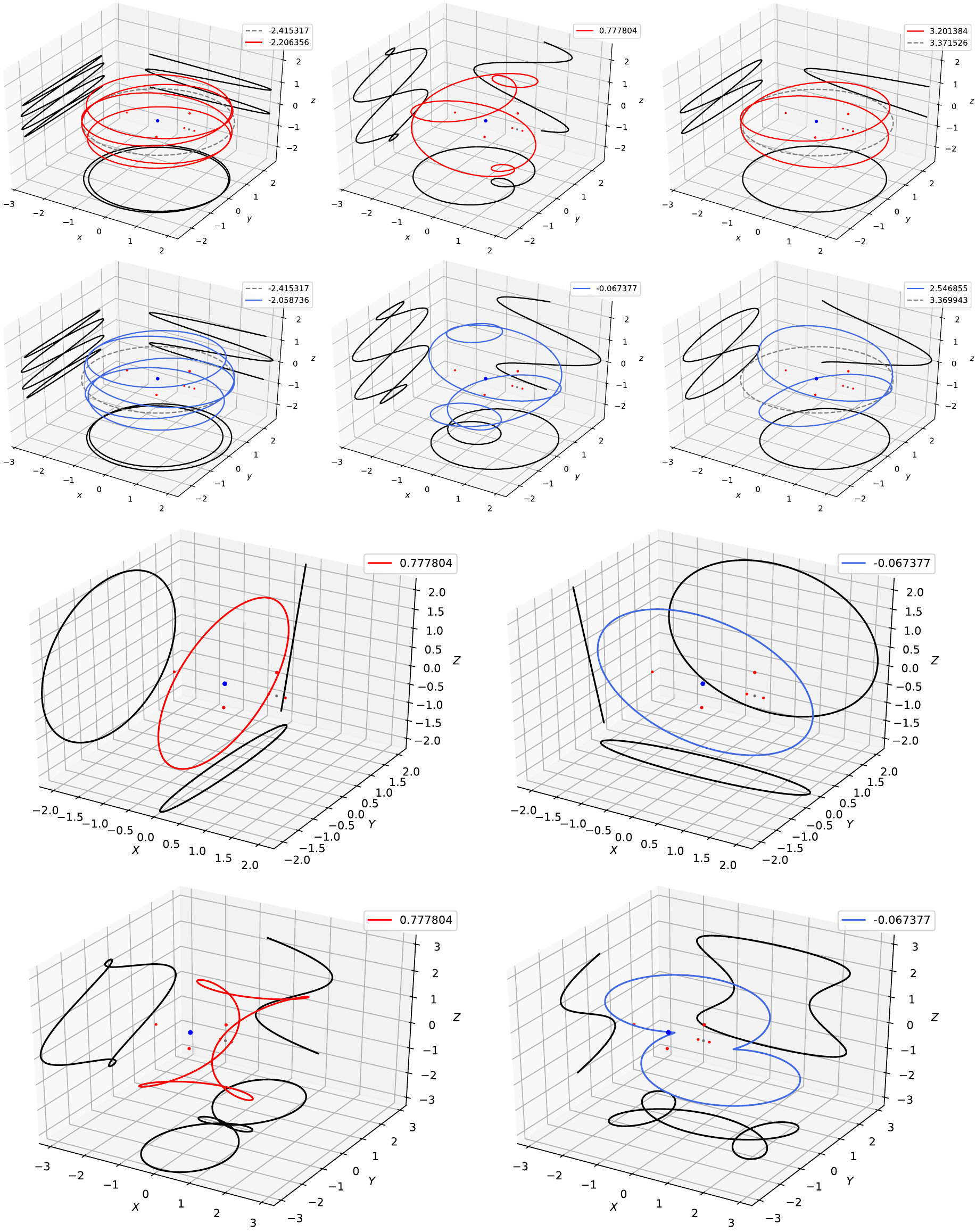}
	\caption{First two rows:\ Orbits of red and blue branch that are doubly symmetric w.r.t.\ the $xz$-plane and the $x$-axis, forming bridge (from left to right) between 1:4 retrograde $\kappa_-$ orbit and 1:2 direct $\kappa_+$ orbit (grey dashed left and right, respectively).\ Orbits in the middle are 1:3 resonant orbits that are plotted in the third row in the inertial frame and in the fourth row in the Moon-centered inertial frame, after three revolutions of the Earth and the Moon.}
	\label{plot_5}
\end{figure}

\begin{figure}[t!]
	\centering
	\includegraphics[width=1\linewidth]{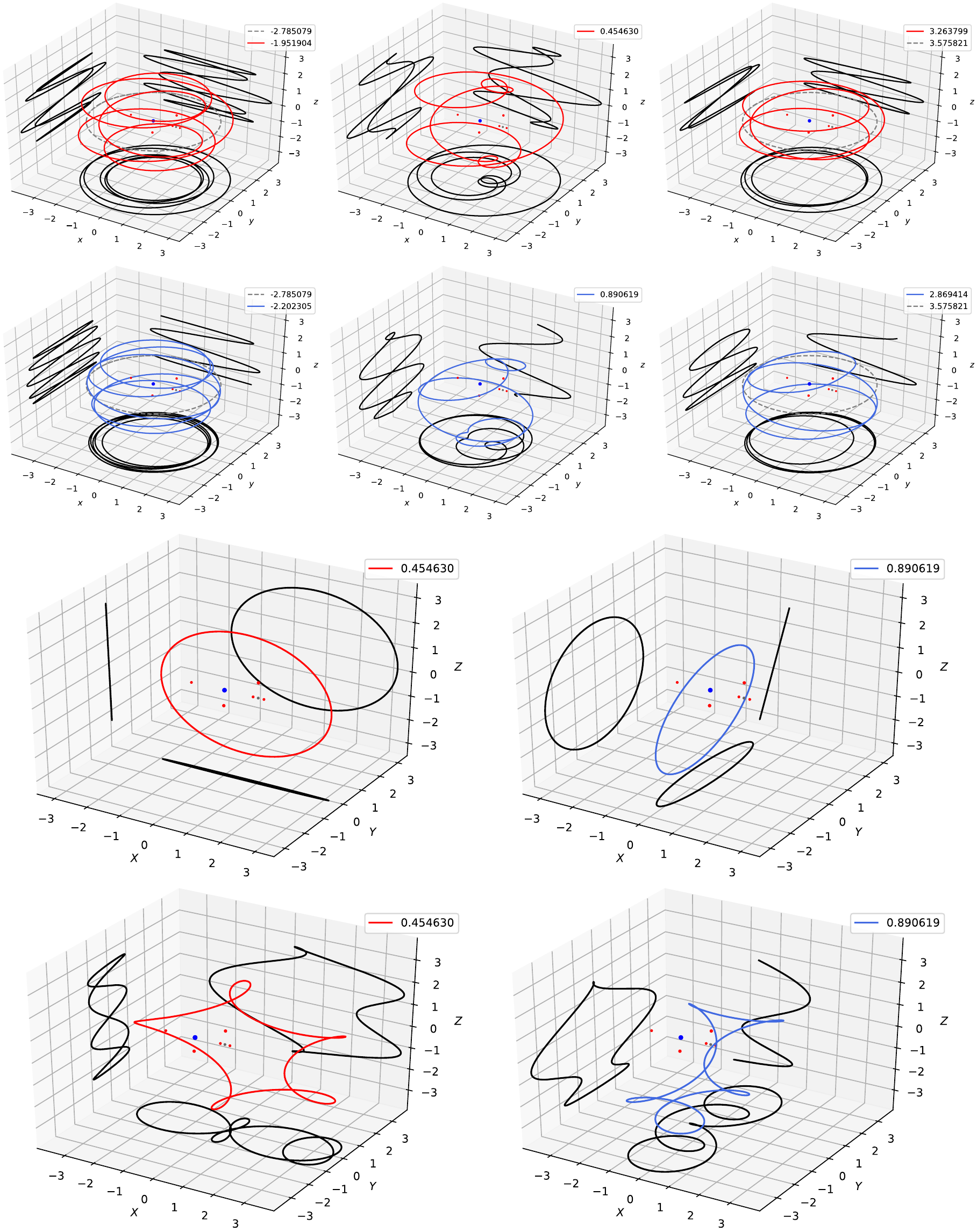}
	\caption{First two rows:\ Orbits of red and blue branch forming bridge, from left to right, between 1:5 retrograde $\kappa_-$ orbit and 1:3 direct $\kappa_+$ orbit (grey dashed left and right, respectively).\ While the red orbits are symmetric w.r.t.\ the $x$-axis, the blue orbits are symmetric w.r.t.\ the $xz$-plane.\ Along the continuation, some loops appear and disappear on the orbits.\ Orbits in the middle are 1:4 resonant orbits that are plotted in the third row in the inertial frame and in the fourth row in the Moon-centered inertial frame, after four revolutions of the Earth and the Moon.}
	\label{plot_6}
\end{figure}

\begin{figure}[t!]
	\centering
	\includegraphics[width=1\linewidth]{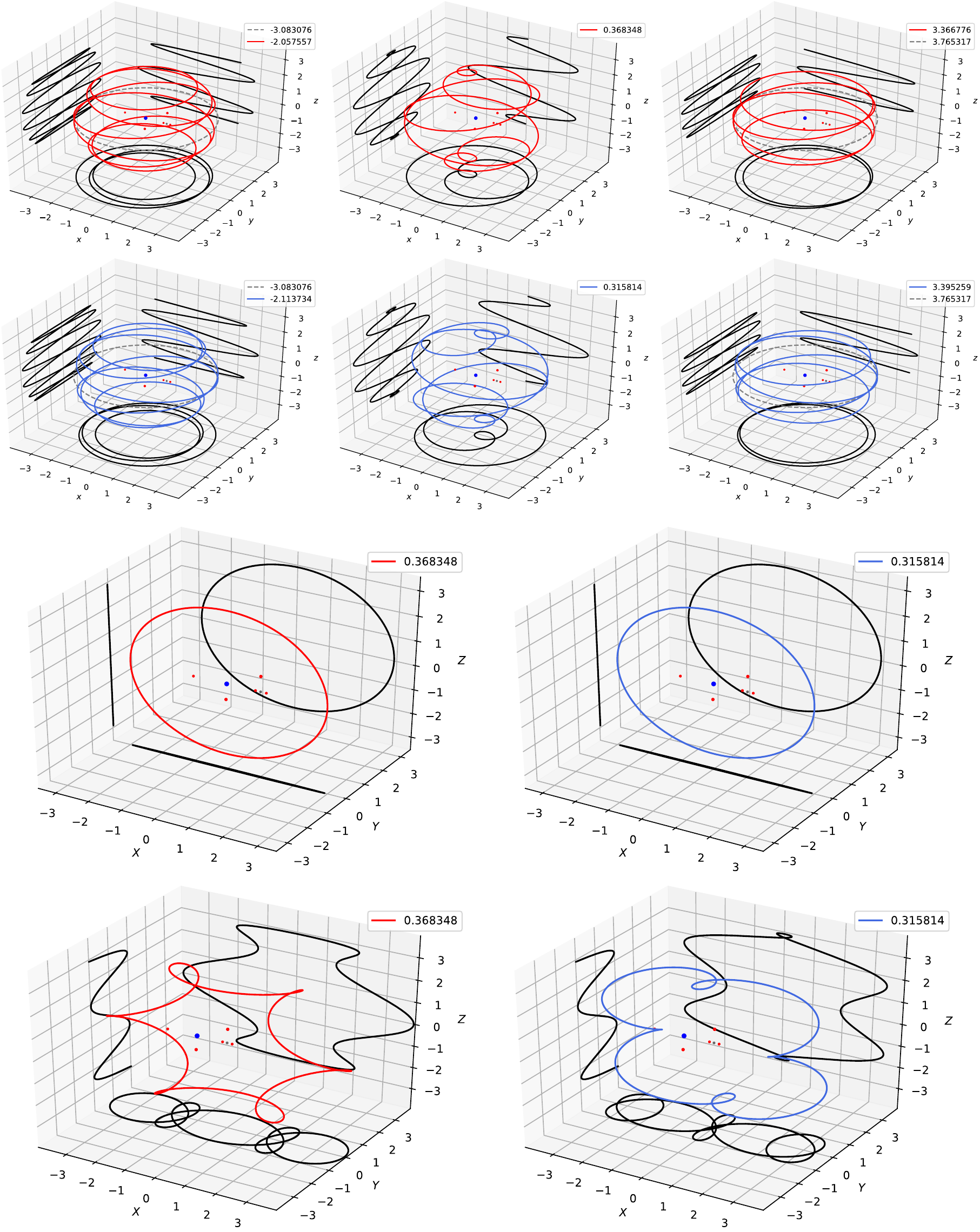}
	\caption{First two rows:\ Orbits of red and blue branch that are doubly symmetric w.r.t.\ the $xz$-plane and the $x$-axis, forming bridge (from left to right) between 1:6 retrograde $\kappa_-$ orbit and 1:4 direct $\kappa_+$ orbit (grey dashed left and right, respectively).\ Along the continuation, some loops appear and disappear on the orbits.\ Orbits in the middle are 1:5 resonant orbits that are plotted in the third row in the inertial frame and in the fourth row in the Moon-centered inertial frame, after five revolutions of the Earth and the Moon.}
	\label{plot_7}
\end{figure}

\noindent
\framebox{\textit{(2) Bifurcation from 1:$k$ retrograde $\kappa_-$ orbits ($k=4,5,6$).}}\ The vertical self-resonant bifurcations of the 1:$k$ retrograde $\kappa_-$ orbits produce all the same results, for $k=4,5,6$.\ For each $k$, the Conley--Zehnder index of the $k$-covering of the 1:$k$ retrograde $\kappa_-$ orbits jumps from 2 to 4, and all the bifurcated branches terminate planar at the $\tilde{k}$-covering of the 1:$\tilde{k}$ direct $\kappa_+$ orbits, where $\tilde{k}=k-2$.\ Therefore, for each $k$, we obtain bridge families between the 1:$k$ retrograde $\kappa_-$ and the 1:$(k-2)$ direct $\kappa_+$ orbits.\ The bifurcation diagrams are shown in Figure \ref{plot_bif_2}, and orbit plots are provided in Figure \ref{plot_5}, Figure \ref{plot_6} and Figure \ref{plot_7}.\ Since, with increasing distance, the periods of the retrograde and direct comet orbits approach $2 \pi$ from below and from above, respectively, every bridge provides 1:$(k-1)$ resonant orbits.\

\noindent
\framebox{\textit{(2a) Bridge between 1:4 retrograde $\kappa_-$ and 1:2 direct $\kappa_+$ orbits.}}\ The period of the 4-covering of the 1:4 retrograde $\kappa_-$ orbit at $C \approx -2.415316$ is $T \approx 18.870232$ ($\approx$ 82.04 days), with the closest approaches of 421382.11 km and 791532.85 km to the lunar and Earth surface, respectively.\ The orbits of the bifurcated red and blue branches are doubly symmetric w.r.t.\ the $xz$-plane and the $x$-axis.\ Some orbits are plotted in Figure \ref{plot_5}, and database are listed in Table \ref{data_4}.\ By inspection the orbit plots in Figure \ref{plot_5}, various loops appear and disappear on the red and blue orbits, and it seems that there is a $yz$-plane symmetry with respect to the origin which relates the red and blue families.\ Along the continuation of the red branch, the orbits with Conley--Zehnder index 2 are stable.\ The vertical initial parameter $z$ reaches a maximum value $z \approx 2.077839$ at $C \approx 0.463212$, then it decreases and terminates planar at the first period-doubling branch point of the 1:2 direct $\kappa_+$ orbit at $C \approx 3.371526$ with period $T \approx 18.817342$ ($\approx$ 81.81 days) and the minimal distances to the lunar and Earth surface of 422072.93 km and 790714.56 km, correspondingly.\ Concerning the blue branch, the initial parameter $\dot{z}$ takes a peak at $\dot{z} \approx 0.693535$ at $C \approx 0.452400$, and then it goes down to plane on the second period-doubling bifurcation point of the 1:2 direct $\kappa_+$ orbit at $C \approx 3.369943$ with period $T \approx 18.840632$ ($\approx$ 81.91 days) and the closest approaches of 420768.32 km and 789384.71 km to the lunar and Earth surface, respectively.\ The members of the blue branch are all unstable (see Table \ref{data_4}).\ Notice that due to the various negative hyperbolic phases of the planar and vertical stability indices (see Figure \ref{plot_2}), the Conley--Zehnder index of the double-covered 1:2 direct $\kappa_+$ orbits first goes from 2 to~$\overline{3}$ (bad orbits), then from $\overline{3}$ to 4, and then from 4 to $\overline{5}$ (bad again), as shown in Figure \ref{plot_bif_2}.\ Furthermore, along the red and blue branches, 1:3 resonant orbits occur that are shown in Figure~\ref{plot_5}.\ While the red 1:3 resonant orbit has Conley--Zehnder index 3 with elliptic and positive hyperbolic Floquet multipliers of the form $\{ e^{\pm i \theta}, \lambda, \lambda^{-1} \}$, with $\theta \approx 5.761$ and $\lambda \approx 1.032$, the blue 1:3 resonant orbit has index 4 with positive hyperbolic Floquet multipliers close to 1 (see Table \ref{data_4}).\
\\
\framebox{\textit{(2b) Bridge between 1:5 retrograde $\kappa_-$ and 1:3 direct $\kappa_+$ orbits.}}\ The 5-covering of the 1:5 retrograde $\kappa_-$ orbit has period $T \approx 25.144807$ ($\approx$ 109.32 days) at $C \approx -2.785078$, with the closest approaches of 589344.48 km and 959656.28 km to the lunar and Earth surface, respectively.\ While the orbits of the red branch are symmetric w.r.t.\ the $x$-axis, the orbits of the blue branch are symmetric w.r.t.\ the $xz$-plane.\ By applying the $\sigma$-symmetry, each of the two branches has a symmetric branch.\ Some orbit plots are provided in Figure \ref{plot_6}, and database can be found in Table \ref{data_5}.\ Both branches terminate at the 3-covering of the 1:3 direct $\kappa_+$ orbit at $C \approx 3.575821$ with period $T \approx 25.120666$ ($\approx$ 109.22 days) and the minimal distances to the lunar and Earth surface of 589000.10 km and 959006.20 km, respectively.\ The blue branch consists entirely of orbits that are of elliptic and positive hyperbolic type, with constant Conley--Zehnder index~3.\ The maximum value of the initial parameter $z$ is $z \approx 1.735593$ at $C \approx -0.817202$.\ Along the continuation of the red branch, the orbits start with index 2, then the index jumps to 3, and then it jumps back to 2.\ The orbits with index 2 are stable.\ The initial parameter $\dot{z}$ reaches a maximum value $\dot{z} \approx 0.812872$ at $C \approx -0.192251$.\ Both branches provide 1:4 resonant orbits (plotted in Figure \ref{plot_6}) that have index 3 with elliptic and positive hyperbolic Floquet multipliers of the form $\{ e^{\pm i \theta}, \lambda, \lambda^{-1} \}$, with $\theta \approx 5.995$ and $\lambda \approx 1.045$ for the red orbit, and with $\theta \approx 5.551$ and $\lambda \approx 1.042$ for the blue orbit.\
\\
\framebox{\textit{(2c) Bridge between 1:6 retrograde $\kappa_-$ and 1:4 direct $\kappa_+$ orbits.}}\ The period of the 6-covering of the 1:6 retrograde $\kappa_-$ orbit at $C \approx -3.083075$ is $T \approx 31.424193$ ($\approx$ 136.63 days) with the minimal distances of 744204.34 km and 1114571.23 km to the lunar and Earth surface, correspondingly.\ The members of the bifurcated red and blue branches are doubly symmetric w.r.t.\ the $xz$-plane and the $x$-axis.\ Some orbits are plotted in Figure \ref{plot_7}, and data are provided in Table \ref{data_6}.\ Both families terminate at the ecliptic on the 4-covering of the 1:4 direct $\kappa_+$ orbit at $C \approx 3.765317$ with period $T \approx 31.407568$ ($\approx$ 136.56 days) and the closest approaches of 743926.32 km and 1114183.43 km to the lunar and Earth surface, respectively.\ Similar as before, in view of the orbit plots in Figure~\ref{plot_7}, it seems that there is a $yz$-plane symmetry with respect to the origin, that connects both branches.\ Orbits of the red branch start with Conley--Zehnder index 2, then jumps to 3, and then jumps back to 2 before reaching the 1:4 direct $\kappa_+$ orbit.\ All the orbits with index 2 are stable, while those with index 3 are of elliptic and positive hyperbolic type.\ The out-of-plane initial parameter $\dot{z}$ reaches a maximum value $\dot{z} \approx 0.585070$ at $C \approx 0.370683$.\ The members of the blue branch start with index 3, then the index jumps to~4, and then jumps back to 3.\ All the orbits are unstable, being elliptic and positive hyperbolic with index 3, and purely positive hyperbolic with index 4.\ The maximum value of the initial parameter $\dot{z}$ is $\dot{z} \approx 0.584669$ at $C \approx 0.356784$.\ Along the continuation of both branches, we find 1:5 resonant orbits that are shown in Figure \ref{plot_7}.\ While the red 1:5 resonant orbit has index 3 with elliptic and positive hyperbolic Floquet multipliers of the form $\{ e^{\pm i \theta}, \lambda, \lambda^{-1} \}$, with $\theta \approx 6.129$ and $\lambda \approx 1.021$, the blue 1:5 resonant orbit has index 4 with all Floquet multipliers in the positive hyperbolic region close to 1 (see Table \ref{data_6}).\

\subsubsection{Single-turn bifurcation from the direct comet orbits}

\begin{figure}[t!]
	\centering
	\includegraphics[width=1\linewidth]{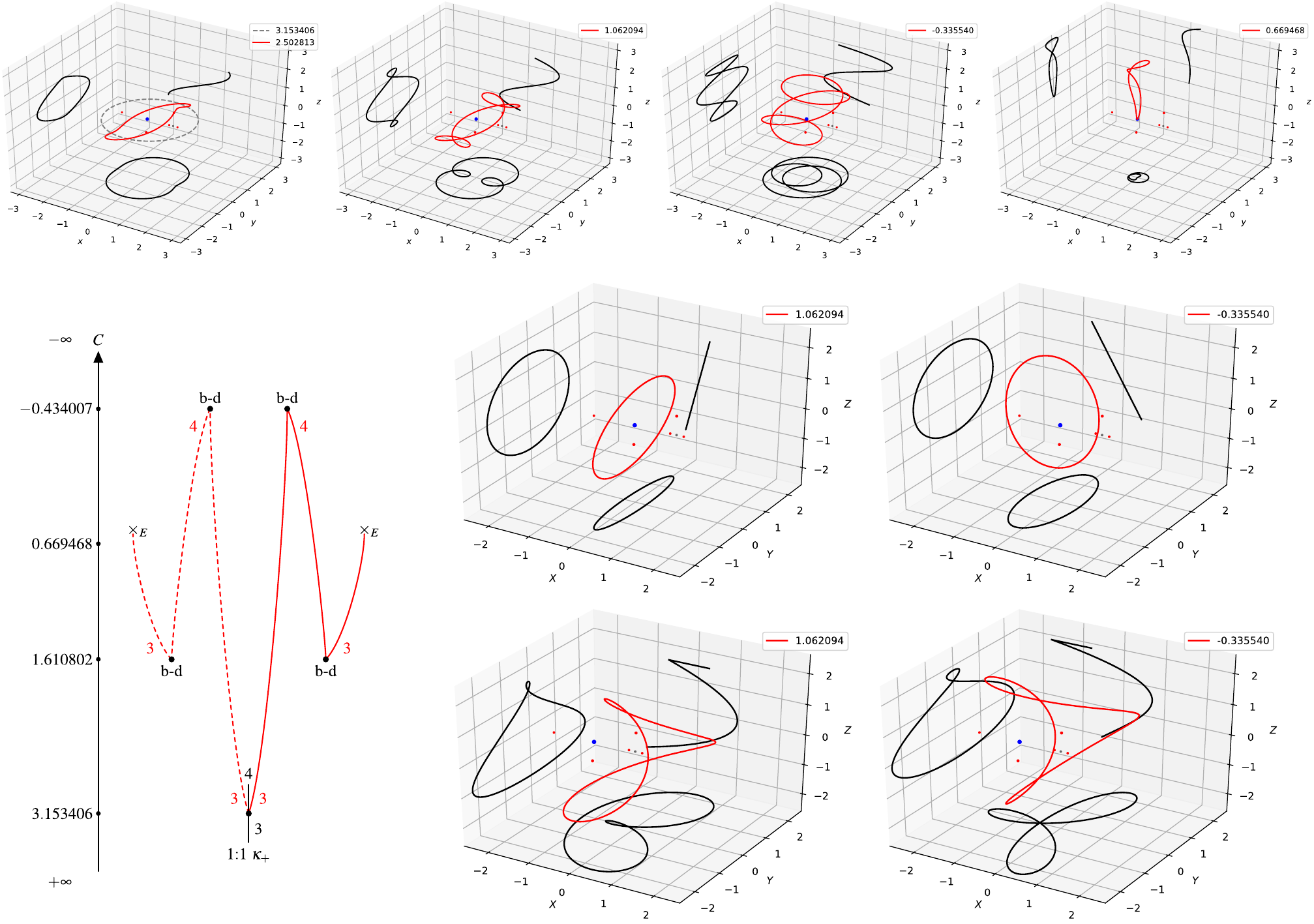}
	\caption{Top, from left to right, shows red orbits (symmetric w.r.t.\ the $xz$-plane) that vertically bifurcate from 1:1 retrograde comet $\kappa_+$ orbit (grey dashed left).\ Left shows corresponding bifurcation diagram that is continued up to approaching collision with the Earth, indicated by ``$\times_E$''.\ Third and fourth row of orbit plots show 1:2 resonant orbits that are plotted at the top in the rotating frame, at the third row in the inertial frame and at the fourth row in the Moon-centered inertial frame, after two revolutions of the Earth and the Moon.}
	\label{plot_8}
\end{figure}

The 1:1 direct $\kappa_+$ orbit at $C \approx 3.153406$ has period $T \approx 12.508392$ ($\approx$ 53.56 days) with minimal distances of 269148.38 km and 572028.09 km to the lunar and Earth surface, respectively.\ At this branch points, we generate a branch whose orbits are simple symmetric w.r.t.\ the $xz$-plane.\ By using the $\sigma$-symmetry, one obtains a symmetric branch.\ The resulting bifurcation graph together with some orbit plots are shown in Figure \ref{plot_8}, and database are provided in Table~\ref{data_7}.\ At the 1:1 direct $\kappa_+$ orbit the Conley--Zehnder index jumps from 3 to 4, and the members of the bifurcated branch start with index 3 being of elliptic and positive hyperbolic type.\ Along the continuation, there are some loops that appear and disappear on the orbits.\ After undergoing two birth-death bifurcations, the orbits approach collision with the Earth.\ At $C \approx 0.669468$ we stopped with the continuation where the orbit has the closest approach of 3523.05 km to the origin of the Earth, which means that it collides with the Earth surface.\ Furthermore, the orbit's geometry in the near-collision phase is of special shape:\ most of the period the orbits are far from the Earth with a large excursion of hook-shaped form in the vertical direction, so that the near-collision moment is very brief, as shown in the right top in Figure \ref{plot_8}.\ This kind of orbit is known for small mass values, as bifurcation from vertical collision orbits in the rotating Kepler problem \cite{belbruno_1,belbruno_2}.\ We also find two 1:2 resonant orbits that are shown in Figure \ref{plot_8}.\ The one at $C \approx 1.062094$ has index 3 with elliptic and positive hyperbolic Floquet multipliers of the form $\{ e^{\pm i \theta}, \lambda, \lambda^{-1} \}$, with $\theta \approx 5.484$ and $\lambda \approx 1.157$.\ The one at $C \approx -0.335539$ has index 4 with positive hyperbolic Floquet multipliers close to 1 (see Table~\ref{data_7}).\

\section{Conclusion}
\label{sec6}

\begin{figure}[t!]
	\centering
	\includegraphics[width=1\linewidth]{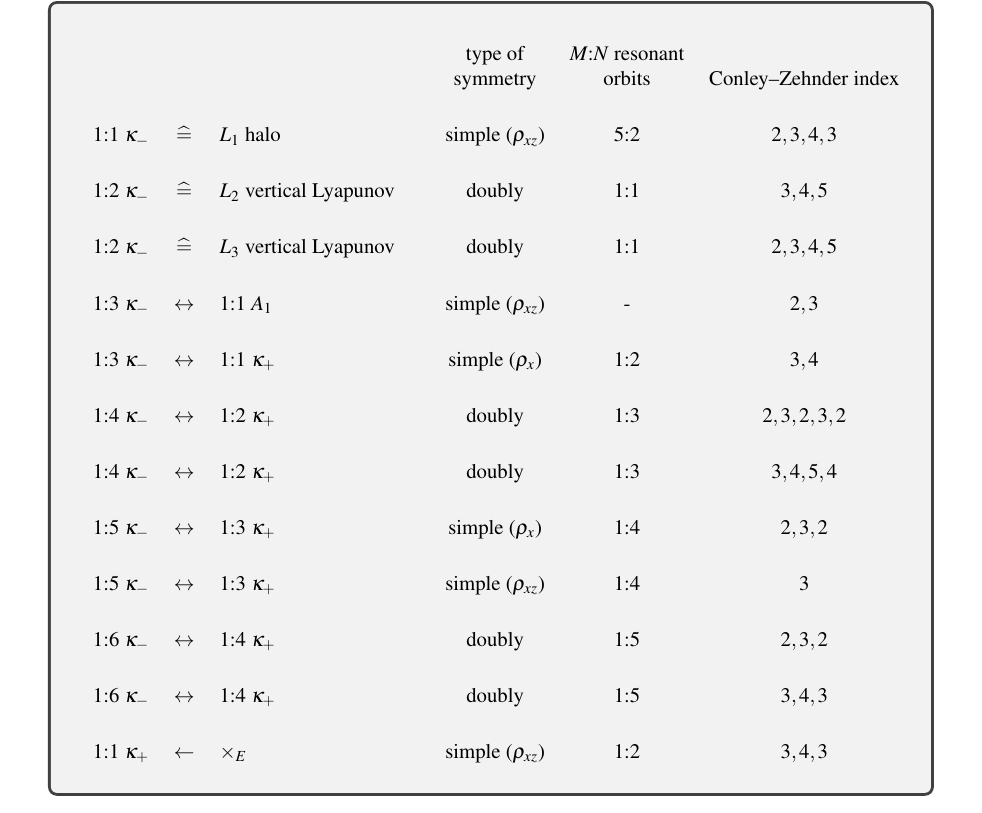}
	\caption{Main bifurcation results.\ ``$\widehat{=}$'' indicates correspondence of the bifurcated branch from the left orbit with the right branch.\ Double arrows indicate bridge families between left and right orbits with corresponding 1:$k$ self-resonance.\ Left arrow indicates non-closed branch that approach collision with the Earth denoted by ``$\times_E$''.\ Families whose orbits have Conley--Zehnder index 2 provide stable solutions.}
	\label{figure_results}
\end{figure}

Within the CR3BP framework, we studied the comet-type periodic orbit families that are generated from very large retrograde and direct Keplerian motions around the common center of mass of the primaries, that we label by $\kappa_{\mp}$, respectively.\ We discussed an analytical approach of their existence in Theorem~\ref{theorem1}, based on the classical Poincaré continuation method, and also determined the Conley--Zehnder index, defined as a Maslov index using a crossing form.\ Then, by applying a standard corrector--predictor algorithm, we explored numerically the two families $\kappa_{\pm}$ and their vertical self-resonant bifurcations in the Earth--Moon environment.\ We computed the stability indices, identified vertical self-resonant orbits up to the multiplicity six, analyzed the vertically bifurcated branches, including bridge families, discussed their orbital characteristics and investigated resonant orbits.\ We computed their Conley--Zehnder indices that provide a systematic classification of the orbit families.\ The main bifurcation results are summarized in Figure \ref{figure_results}, and detailed descriptions are provided in Section \ref{sec5}.\ We notice that the bifurcation structure of the bridge families among the retrograde and direct comet orbits is similar to that of the retrograde and direct orbits that are generated by small Keplerian motions near the primaries \cite{aydin_babylonian,aydin_cz,aydin_batkhin,aydin_dro,moreno_aydin}, namely two more coverings of the retrograde are necessary to form a bridge to the direct orbits.\ In contrast to the small retrograde and direct orbits, where the Conley--Zehnder index of the multiple-coverings becomes very high, the multiple-coverings of the comet orbits have an index jump from 2 to 4 in all considered~cases.\

\clearpage
\begin{appendix}
\setcounter{secnumdepth}{0}
\renewcommand{\thetable}{A\arabic{table}}
\section{Appendix Tables of data}
%\label{appendix}
	
The database contains initial data, Floquet multipliers with corresponding Krein signatures and Conley--Zehnder indices of selected orbits associated to the families we have studied.\ For a pair of hyperbolic Floquet multipliers we denote by $\lambda$ the one with $|\lambda| > 1$.\ For a pair of elliptic Floquet multipliers we denote by~$\theta$ the rotation angle modulo $2 \pi$ characterized by Krein signature.\ Horizontal dashed lines indicate birth-death bifurcations.\
	
\begin{table}[ht]\fontsize{10}{10}\selectfont \centering
		\caption{Data for $\kappa_-$ (first block) and $\kappa_+$ (second block).}
		\begin{tabular}{c|c|c|c|c|c}
			$C$ & $x(0)$ & $\dot{y}(0)$ & $T/2$ & Floquet multipliers \& Krein sign & $\mu_{CZ}$ / $\mu_{CZ}^p$ / $\mu_{CZ}^s$\\
			\hline $-$3.73079020 & 3.96375030 & $-$4.46622787 & 2.788167 & $(+)$ $\theta_p \approx 0.706$, $(+)$ $\theta_s \approx 0.707$ & 2 / 1 / 1\\
			$-$3.08307598 & 2.92838023 & $-$3.51324570 & 2.618682 & $(+)$ $\theta_p \approx 1.044$, $(+)$ $\theta_s \approx 2 \pi/6$ & 2 / 1 / 1\\
			$-$2.78507825 & 2.52551891 & $-$3.15557876 & 2.514480 & $(+)$ $\theta_p \approx 1.251$, $(+)$ $\theta_s \approx 2 \pi/5$ & 2 / 1 / 1\\
			$-$2.41531689 & 2.08857210 & $-$2.78208822 & 2.358779 & $(+)$ $\theta_p \approx 1.560$, $(+)$ $\theta_s \approx 2 \pi/4$ & 2 / 1 / 1\\
			$-$1.91761341 & 1.60566275 & $-$2.39923349 & 2.103212 & $(+)$ $\theta_p \approx 2.059$, $(+)$ $\theta_s \approx 2 \pi/3$ & 2 / 1 / 1\\
			$-$1.33311990 & 1.17816938 & $-$2.12337566 & 1.742922 & $(+)$ $\theta_p \approx 2.650$, $\lambda_s \approx -1.000$ & 2 / 1 / 1\\
			$-$1.19019267 & 1.09454322 & $-$2.09792012 & 1.634914 & $(+)$ $\theta_p \approx 2.741$, $\lambda_s \approx -1.397$ & 2 / 1 / 1\\
			$-$1.12676203 & 1.06081793 & $-$2.10393894 & 1.574541 & $(+)$ $\theta_p \approx 2.934$, $\lambda_s \approx -1.000$ & 2 / 1 / 1\\
			$-$1.06179859 & 1.03118167 & $-$2.14002812 & 1.490894 & $\lambda_p \approx -2.251$, $(-)$ $\theta_s \approx 2 \pi$ $2/3$ & 2 / 1 / 1\\
			$-$1.04060107 & 1.02363545 & $-$2.16219513 & 1.456837 & $\lambda_p \approx -3.079$, $(-)$ $\theta_s \approx 2 \pi$ $3/4$ & 2 / 1 / 1\\
			$-$1.01611986 & 1.01670394 & $-$2.19373120 & 1.414661 & $\lambda_p \approx -4.550$, $\lambda_s \approx 1.000$ & $2 \shortto 3$ / 1 / $1 \shortto 2$\\
			$-$0.48769543 & $-$0.27646347 & 2.83871407 & 0.796606 & $\lambda_p \approx -57.70$, $\lambda_s \approx 300.5$ & 3 / 1 / 2\\
			$-$0.41303508 & $-$0.17803693 & 3.51792654 & 0.671272 & $(+)$ $\theta_p \approx 2.760$, $\lambda_s \approx 533.5$ & 3 / 1 / 2\\
			$-$0.41292572 & $-$0.17508887 & 3.54826837 & 0.667245 & $(+)$ $\theta_p \approx 0.019$, $\lambda_s \approx 544.3$ & 3 / 1 / 2\\
			\hdashline $-$0.41292568 & $-$0.17508886 & 3.54826847 & 0.667244 & $\lambda_p \approx 1.004$, $\lambda_s \approx 544.3$ & 2 / 0 / 2\\
			$-$0.95924764 & $-$0.06168512 & 6.39335159 & 0.477327 & $\lambda_p \approx 848.9$, $\lambda_s \approx 1926$ & 2 / 0 / 2\\
			$-$1.68063393 & $-$0.04012751 & 8.50504725 & 0.417761 & $\lambda_p \approx 2135$, $\lambda_s \approx 3554$ & 2 / 0 / 2\\
			$-$2.51369826 & $-$0.02999931 & 10.64223834 & 0.378057 & $\lambda_p \approx 4033$, $\lambda_s \approx 5797$ & 2 / 0 / 2\\
			\hline 3.94727405 & 3.32100001 & $-$2.77216975 & 3.763832 & $(+)$ $\theta_p \approx 1.243$, $(+)$ $\theta_s \approx 1.245$ & 2 / 1 / 1\\
			3.76531791 & 2.92763105 & $-$2.34313685 & 3.925946 & $(+)$ $\theta_p \approx 1.566$, $(+)$ $\theta_s \approx 2 \pi/4$ & 2 / 1 / 1\\
			3.57582149 & 2.52462303 & $-$1.89539901 & 4.186777 & $(+)$ $\theta_p \approx 2.086$, $(+)$ $\theta_s \approx 2 \pi/3$ & 2 / 1 / 1\\
			3.37152656 & 2.09036920 & $-$1.39994341 & 4.704335 & $(+)$ $\theta_p \approx 3.134$, $\lambda_s \approx -1.000$ & 2 / 1 / 1\\
			3.37139816 & 2.09009402 & $-$1.39962427 & 4.704806 & $\lambda_p \approx -1.000$, $\lambda_s \approx -1.003$ & 2 / 1 / 1\\
			3.36994316 & 2.08697533 & $-$1.39600678 & 4.710158 & $\lambda_p \approx -1.021$, $\lambda_s \approx -1.000$ & 2 / 1 / 1\\
			3.36445490 & 2.07520682 & $-$1.38234813 & 4.730661 & $\lambda_p \approx -1.000$, $(-)$ $\theta_s \approx 3.188$ & 2 / 1 / 1\\
			3.26371835 & 1.85842481 & $-$1.12877794 & 5.226019 & $(-)$ $\theta_p \approx 4.150$, $(-)$ $\theta_s \approx 2 \pi$ $2/3$ & 2 / 1 / 1\\
			3.22782403 & 1.78256000 & $-$1.03977090 & 5.484844 & $(-)$ $\theta_p \approx 4.668$, $(-)$ $\theta_s \approx 2 \pi$ $3/4$ & 2 / 1 / 1\\
			3.15340606 & 1.69254257 & $-$0.95118846 & 6.254196 & $\lambda_p \approx 1.473$, $\lambda_s \approx 1.000$ & $3 \shortto 4$ / 2 / $1 \shortto 2$\\
			3.15194964 & 1.69613732 & $-$0.95696577 & 6.268160 & $\lambda_p \approx 1.498$, $\lambda_s \approx 1.017$ & 4 / 2 / 2\\
			3.15136315 & 1.69775817 & $-$0.95952922 & 6.273520 & $\lambda_p \approx 1.508$, $\lambda_s \approx 1.000$ & $4 \shortto 5$ / 2 / $2 \shortto 3$\\
			3.15055237 & 1.70016298 & $-$0.96329648 & 6.280659 & $\lambda_p \approx 1.521$, $(+)$ $\theta_s \approx 0.027$ & 5 / 2 / 3\\
			2.98165907 & 2.17197153 & $-$1.63122614 & 6.324042 & $\lambda_p \approx 3.678$, $(+)$ $\theta_s \approx 0.106$ & 5 / 2 / 3\\
			1.49947321 & 3.07313502 & $-$2.93201321 & 6.274403 & $\lambda_p \approx 2.319$, $(+)$ $\theta_s \approx 0.016$ & 5 / 2 / 3
		\end{tabular}
		\label{data_1}
\end{table}

\begin{table}[ht]\fontsize{10}{10}\selectfont \centering
		\caption{Data for orbit families branching off from 1:1 $\kappa_-$that correspond to $L_1$ halo (first block), and from~1:2 $\kappa_-$ that correspond to $L_2$ and $L_3$ vertical Lyapunov (second and third block, respectively).}
		\begin{tabular}{c|c|c|c|c|c|c}
			$C$ & $x(0)$ & $z(0)$ & $\dot{y}(0)$ & $T/2$ & Floquet multipliers \& Krein sign & $\mu_{CZ}$\\
			\hline $-$1.01611986 & $-$0.8457682 & 0 & 2.02848068 & 1.414661 & $\lambda_p \approx -4.550$, $\lambda_s \approx 1.000$ & $2 \shortto 3$\\
			$-$1.01121740 & $-$0.8480226 & 0.04500000 & 2.02578862 & 1.417869 & $\lambda \approx -4.980$, $(-)$ $\theta \approx 5.264$ & 2\\
			$-$0.98614733 & $-$0.8605044 & 0.11500000 & 2.01183966 & 1.435595 & $\lambda_1 \approx -7.801$, $\lambda_2 \approx -1.194$ & 2\\
			1.00758977 & $-$0.0058898 & 0.99501230 & 0.99764024 & 1.556055 & $513.3 \pm 376.2 i$, $0.001 \pm 0.001 i$ & 2\\
			2.85670417 & 0.90106238 & 0.38527217 & 0.10036771 & $\approx 2 \pi / 5$ & $12.24 \pm 16.43 i$, $0.029 \pm 0.039 i$ & 2\\	
			2.94065356 & 0.92852311 & 0.29459187 & 0.08181840 & 1.083121 & $(-)$ $\theta_1 \approx 4.661$, $(+)$ $\theta_2 \approx 1.626$ & 2\\
			2.94342187 & 0.92915408 & 0.29140652 & 0.08174723 & 1.075469 & $(-)$ $\theta \approx 5.504$, $\lambda \approx -1.000$ & 2\\
			2.94699116 & 0.92992363 & 0.28728227 & 0.08174985 & 1.065406 & $\lambda_1 \approx 1.000$, $\lambda_2 \approx -2.623$ & $2 \shortto 3$\\
			2.97165710 & 0.93338449 & 0.25816317 & 0.08596446 & 0.990270 & $\lambda_1 \approx 2.373$, $\lambda_2 \approx -5.415$ & 3\\
			3.00400727 & 0.91276840 & 0.20718952 & 0.15444698 & 0.915881 & $\lambda_1 \approx 1.001$, $\lambda_2 \approx -4.433$ & 3\\
			\hdashline 3.00400727 & 0.91276832 & 0.20718947 & 0.15444717 & 0.915881 & $(+)$ $\theta \approx 0.001$, $\lambda \approx -4.433$ & 4\\
			2.99862220 & 0.88116552 & 0.19373322 & 0.22165687 & 1.056634 & $(+)$ $\theta \approx 0.924$, $\lambda \approx -1.000$ & 4\\
			2.99783519 & 0.87193582 & 0.19019461 & 0.23730262 & 1.115036 & $(+)$ $\theta_1 \approx 0.002$, $(+)$ $\theta_2 \approx 1.957$ & 4\\
			\hdashline 2.99783519 & 0.87193578 & 0.19019460 & 0.23730267 & 1.115037 & $\lambda \approx 1.002$, $(+)$ $\theta \approx 1.957$ & 3\\
			3.00495876 & 0.85258364 & 0.17826198 & 0.26125187 & $\approx 2 \pi / 5$ & $\lambda \approx 17.16$, $(+)$ $\theta \approx 2.438$ & 3\\
			3.11960141 & 0.82631121 & 0.08588652 & 0.20029568 & 1.389850 & $\lambda \approx 838.4$, $(-)$ $\theta \approx 5.342$ & 3\\
			3.17439048 & 0.82336756 & 0 & 0.12634016 & 1.371476 & $\lambda_p \approx 2361$, $\lambda_s \approx 1.000$ & $3 \shortto 4$\\
			\hline $C$ & $x(0)$ & $\dot{y}(0)$ & $\dot{z}(0)$ & $T/4$ & Floquet multipliers \& Krein sign & $\mu_{CZ}$\\
			\hline $-$1.12676203 & 1.06081793 & $-$2.1039389 & 0 & 1.574541 & $(-)$ $\theta_p \approx 5.868$, $\lambda_s \approx 1.000$ & $3 \shortto 4$\\
			$-$1.12353182 & 1.06075604 & $-$2.1025131 & 0.05499999 & 1.574389 & $(-)$ $\theta \approx 5.942$, $\lambda \approx 1.406$ & 3\\
			$-$1.11699017 & 1.06077019 & $-$2.0994238 & 0.09698150 & 1.574377 & $\lambda_1 \approx 1.000$, $\lambda_2 \approx 1.819$ & $3 \shortto 4$\\
			0.16768935 & 1.06471626 & $-$1.4937689 & 0.94100382 & $\approx 2 \pi / 4$ & $\lambda_1 \approx 13.87$, $\lambda_2 \approx 112.7$ & 4\\
			1.00449875 & 1.06930182 & $-$1.1017689 & 1.02485766 & 1.565515 & $\lambda_1 \approx 17.71$, $\lambda_2 \approx 245.3$ & 4\\
			2.96706753 & 1.11194442 & $-$0.1811822 & 0.43588315 & 1.105549 & $\lambda_1 \approx 1.000$, $\lambda_2 \approx 335.3$ & $4 \shortto 5$\\
			3.07784218 & 1.13162660 & $-$0.0676489 & 0.30751995 & 0.935413 & $(+)$ $\theta \approx 0.613$, $\lambda \approx 823.9$ & 5\\
			%3.14423876 & 1.14839882 & $-$0.0173390 & 0.16751995 & 0.890856 & $(+)$ $\theta \approx 0.401$, $\lambda \approx 1550$ & 5\\
			3.16764037 & 1.15450857 & $-$0.0026983 & 0.06751995 & 0.881078 & $(+)$ $\theta \approx 0.294$, $\lambda \approx 1908$ & 5\\
			\hline $C$ & $x(0)$ & $z(0)$ & $\dot{y}(0)$ & $T/4$ & Floquet multipliers \& Krein sign & $\mu_{CZ}$\\
			\hline $-$1.33311990 & 1.17816938 & 0 & $-$2.1233756 & 1.742922 & $(-)$ $\theta_p \approx 5.300$, $\lambda_s \approx 1.000$ & $2 \shortto 3$\\
			$-$1.29018584 & 1.15485565 & 0.10500000 & $-$2.1055272 & 1.729236 & $(-)$ $\theta_1 \approx 5.494$, $(-)$ $\theta_2 \approx 5.725$ & 2\\
			$-$1.08659142 & 1.04597982 & 0.24798430 & $-$2.0233626 & 1.658151 & $\lambda \approx 1.000$, $(-)$ $\theta \approx 5.694$ & $2 \shortto 3$\\
			$-$0.88712597 & 0.94352218 & 0.38500000 & $-$1.9384665 & 1.601378 & $\lambda \approx 1.121$, $(-)$ $\theta \approx 5.970$ & 3\\
			$-$0.64402169 & 0.82132984 & 0.56898610 & $-$1.8213538 & 1.577771 & $\lambda_1 \approx 1.577$, $\lambda_2 \approx 1.000$ & $3 \shortto 4$\\
			$-$0.39249360 & 0.69516788 & 0.71121000 & $-$1.6961804 & $\approx 2 \pi / 4$ & $\lambda_1 \approx 1.333$, $\lambda_2 \approx 1.039$ & 4\\
			0.03075099 & 0.48299074 & 0.86792609 & $-$1.4845700 & 1.566650 & $\lambda_1 \approx 1.605$, $\lambda_2 \approx 1.000$ & $4 \shortto 5$\\
			1.00587257 & $-$0.00506155 & 0.99597899 & $-$0.9975243 & 1.563844 & $\lambda \approx 2.106$, $(+)$ $\theta \approx 0.029$ & 5\\
			2.62321333 & $-$0.81165561 & 0.58962477 & $-$0.1928997 & 1.562603 & $\lambda \approx 2.861$, $(+)$ $\theta \approx 0.032$ & 5\\
			%2.01924928 & $-$0.51088663 & 0.86011032 & $-$0.4928997 & 1.562901 & $\lambda \approx 2.582$, $(+)$ $\theta \approx 0.032$ & 5\\
			%2.83771288 & $-$0.91834774 & 0.40635136 & $-$0.0864872 & 1.562522 & $\lambda \approx 2.959$, $(+)$ $\theta \approx 0.032$ & 5\\
			3.00893966 & $-$1.00346817 & 0.05635136 & $-$0.0015921 & 1.562464 & $\lambda \approx 3.038$, $(+)$ $\theta \approx 0.031$ & 5
		\end{tabular}
		\label{data_2}
\end{table}
	
\begin{table}[ht]\fontsize{10}{10}\selectfont \centering
		\caption{Data for two families emerging from 1:3 $\kappa_-$ (orbits in the first block terminate at 1:1 $A_1$, and orbits in the second block terminate at 1:1 $\kappa_+$).}
		\begin{tabular}{c|c|c|c|c|c|c}
			$C$ & $x(0)$ & $z(0)$ & $\dot{y}(0)$ & $T/2$ & Floquet multipliers \& Krein sign & $\mu_{CZ}$\\
			\hline $-$1.91761341 & 1.60566275 & 0 & $-$2.3992334 & 6.309637 & $(-)$ $\theta_p \approx 6.177$, $\lambda_s \approx 1.000$ & $2 \shortto 4$\\
			$-$1.81735436 & 1.45213275 & 0.40000000 & $-$2.2950540 & 6.307873 & $(-)$ $\theta_1 \approx 6.120$, $(-)$ $\theta_2 \approx 6.266$ & 2\\
			$-$1.64473200 & 1.19211399 & 0.53509332 & $-$2.1465201 & 6.305408 & $(-)$ $\theta_1 \approx 5.968$, $(-)$ $\theta_2 \approx 6.257$ & 2\\
			$-$1.55297424 & 0.89906119 & 0.32676374 & $-$2.1151854 & 6.322365 & $(-)$ $\theta_1 \approx 5.964$, $(-)$ $\theta_2 \approx 6.280$ & 2\\
			\hdashline $-$1.55297424 & 0.89697817 & 0.32348360 & $-$2.1159650 & 6.322815 & $(-)$ $\theta \approx 5.689$, $\lambda \approx 1.002$ & 3\\
			$-$1.56235471 & 0.81383824 & 0.13914529 & $-$2.1662091 & 6.358972 & $(-)$ $\theta \approx 5.470$, $\lambda \approx 1.034$ & 3\\
			$-$1.56849475 & 0.79468554 & 0 & $-$2.1850742 & 6.379428 & $(-)$ $\theta_p \approx 5.393$, $\lambda_s \approx 1.000$ & $2 \shortto 3$\\
			\hline $C$ & $x(0)$ & $\dot{y}(0)$ & $\dot{z}(0)$ & $T/2$ & Floquet multipliers \& Krein sign & $\mu_{CZ}$\\
			\hline $-$1.91761341 & 1.60566275 & $-$2.3992334 & 0 & 6.309637 & $(-)$ $\theta_p \approx 6.177$, $\lambda_s \approx 1.000$ & $2 \shortto 4$\\
			$-$1.81392539 & 1.70156027 & $-$2.4199511 & 0.20000000 & 6.307532 & $(-)$ $\theta \approx 6.118$, $\lambda \approx 1.017$ & 3\\
			$-$1.30597461 & 2.02602529 & $-$2.5044415 & 0.36232839 & 6.293265 & $\lambda_1 \approx 1.000$, $\lambda_2 \approx 1.104$ & $3 \shortto 4$\\
			0.01978289 & 1.66817421 & $-$1.8535998 & 0.73396077 & 6.283455 & $\lambda_1 \approx 1.080$, $\lambda_2 \approx 1.772$ & 4\\
			0.19787704 & 1.63039953 & $-$1.7657998 & 0.76349562 & $2 \pi$ & $\lambda_1 \approx 1.091$, $\lambda_2 \approx 1.898$ & 4\\
			0.84557994 & 1.55885049 & $-$1.4935998 & 0.80855708 & 6.281744 & $\lambda_1 \approx 1.131$, $\lambda_2 \approx 2.425$ & 4\\
			2.33904193 & 1.60479821 & $-$1.0768130 & 0.58143972 & 6.276518 & $\lambda_1 \approx 1.282$, $\lambda_2 \approx 2.617$ & 4\\
			2.94488745 & 1.67186668 & $-$0.9839320 & 0.30143972 & 6.274243 & $\lambda_1 \approx 1.374$, $\lambda_2 \approx 1.807$ & 4\\
			3.15136315 & 1.69775817 & $-$0.9595292 & 0 & 6.273520 & $\lambda_p \approx 1.508$, $\lambda_s \approx 1.000$ & $5 \shortto 4$
		\end{tabular}
		\label{data_3}
\end{table}

\begin{table}[ht]\fontsize{10}{10}\selectfont \centering
	\caption{Data for two families bifurcating from 1:4 $\kappa_-$ (both families terminate at 1:2 $\kappa_+$).}
	\begin{tabular}{c|c|c|c|c|c|c}
		$C$ & $x(0)$ & $z(0)$ & $\dot{y}(0)$ & $T/4$ & Floquet multipliers \& Krein sign & $\mu_{CZ}$\\
		\hline $-$2.41531689 & $-$2.0878675 & 0 & 2.78104740 & 4.717558 & $(-)$ $\theta_p \approx 6.241$, $\lambda_s \approx 1.000$ & $2 \shortto 4$\\
		$-$2.20635533 & $-$1.9383856 & 0.77200000 & 2.63130893 & 4.716830 & $(-)$ $\theta_1 \approx 6.194$, $(-)$ $\theta_2 \approx 6.256$ & 2\\
		$-$1.89791187 & $-$1.7169408 & 1.18249000 & 2.40962414 & 4.715987 & $(-)$ $\theta \approx 6.093$, $\lambda \approx 1.000$ & $2 \shortto 3$\\
		$-$1.55961509 & $-$1.4733833 & 1.47271735 & 2.16593152 & 4.715271 & $(-)$ $\theta \approx 6.011$, $\lambda \approx 1.013$ & 3\\
		$-$1.01534919 & $-$1.0807711 & 1.77861858 & 1.77327052 & 4.714397 & $(-)$ $\theta \approx 5.915$, $\lambda \approx 1.000$ & $3 \shortto 2$\\
		$-$0.57577415 & $-$0.7633845 & 1.93482010 & 1.45593152 & 4.713838 & $(-)$ $\theta_1 \approx 5.860$, $(-)$ $\theta_2 \approx 6.273$ & 2\\
		$-$0.16974581 & $-$0.4701675 & 2.02514833 & 1.16279452 & 4.713384 & $(-)$ $\theta \approx 5.821$, $\lambda \approx 1.000$ & $2 \shortto 3$\\
		0.46321270 & $-$0.0131566 & 2.07783917 & 0.70593152 & 4.712724 & $(-)$ $\theta \approx 5.776$, $\lambda \approx 1.024$ & 3\\
		0.77780422 & 0.21391077 & 2.06645864 & 0.47893152 & $\approx 2 \pi$ $3/4$ & $(-)$ $\theta \approx 5.761$, $\lambda \approx 1.032$ & 3\\
		1.23938644 & 0.54697989 & 2.00390183 & 0.14593152 & 4.711854 & $(-)$ $\theta \approx 5.744$, $\lambda \approx 1.042$ & 3\\
		2.82478178 & 1.69181204 & 1.21295943 & $-$0.9996825 & 4.708275 & $(-)$ $\theta \approx 5.822$, $\lambda \approx 1.000$ & $3 \shortto 2$\\
		3.20138458 & 1.96566673 & 0.70006443 & $-$1.2744995 & 4.705982 & $(-)$ $\theta_1 \approx 5.973$, $(-)$ $\theta_2 \approx 6.243$ & 2\\
		3.37152656 & 2.09036920 & 0 & $-$1.3999434 & 4.704335 & $(-)$ $\theta_p \approx 6.269$, $\lambda_s \approx 1.000$ & $2 \shortto \overline{3}$\\
		\hline $C$ & $x(0)$ & $\dot{y}(0)$ & $\dot{z}(0)$ & $T/4$ & Floquet multipliers \& Krein sign & $\mu_{CZ}$\\
		\hline $-$2.41531689 & 2.08857210 & $-$2.7820882 & 0 & 4.717558 & $(-)$ $\theta_p \approx 6.241$, $\lambda_s \approx 1.000$ & $2 \shortto 4$\\
		$-$2.23032078 & 2.08707669 & $-$2.7369527 & 0.24199999 & 4.716389 & $(-)$ $\theta \approx 6.258$, $\lambda \approx 1.082$ & 3\\
		$-$2.05873588 & 2.08606393 & $-$2.6952472 & 0.33133100 & 4.715605 & $\lambda_1 \approx 1.000$, $\lambda_2 \approx 1.152$ & $3 \shortto 4$\\
		$-$1.31916187 & 2.08362347 & $-$2.5163007 & 0.54200000 & 4.713737 & $\lambda_1 \approx 1.028$, $\lambda_2 \approx 1.408$ & 4\\
		$-$0.06737689 & 2.08217994 & $-$2.2146922 & 0.68084679 & $\approx 2 \pi$ $3/4$ & $\lambda_1 \approx 1.008$, $\lambda_2 \approx 1.741$ & 4\\
		0.02928745 & 2.08214076 & $-$2.1914482 & 0.68495225 & 4.712318 & $\lambda_1 \approx 1.000$, $\lambda_2 \approx 1.762$ & $4 \shortto 5$\\
		0.45240023 & 2.08206794 & $-$2.0897922 & 0.69353564 & 4.712029 & $(+)$ $\theta \approx 0.016$, $\lambda \approx 1.846$ & 5\\
		1.87581185 & 2.08305464 & $-$1.7492002 & 0.60725253 & 4.711176 & $\lambda_1 \approx 1.000$, $\lambda_2 \approx 1.952$ & $5 \shortto 4$\\
		2.54685508 & 2.08433874 & $-$1.5897922 & 0.48448151 & 4.710758 & $\lambda_1 \approx 1.035$, $\lambda_2 \approx 1.820$ & 4\\
		3.36994316 & 2.08697533 & $-$1.3960067 & 0 & 4.710158 & $\lambda_p \approx 1.042$, $\lambda_s \approx 1.000$ & $4 \shortto \overline{5}$
	\end{tabular}
	\label{data_4}
\end{table}

\begin{table}[ht]\fontsize{10}{10}\selectfont \centering
	\caption{Data for two families branching off from 1:5 $\kappa_-$ (both families terminate at 1:3 $\kappa_+$).}
	\begin{tabular}{c|c|c|c|c|c|c}
		$C$ & $x(0)$ & $\dot{y}(0)$ & $\dot{z}(0)$ & $T/2$ & Floquet multipliers \& Krein sign & $\mu_{CZ}$\\
		\hline $-$2.78507825 & $-$2.5252333 & 3.15514897 & 0 & 12.572403 & $(-)$ $\theta_p \approx 6.259$, $\lambda_s \approx 1.000$ & $2 \shortto 4$\\
		$-$1.95190358 & $-$2.0608515 & 2.63081580 & 0.50015919 & 12.570616 & $(-)$ $\theta_1 \approx 6.039$, $(-)$ $\theta_2 \approx 6.278$ & 2\\
		$-$1.64857749 & $-$1.6748075 & 2.28599260 & 0.65207487 & 12.570488 & $(-)$ $\theta \approx 5.994$, $\lambda \approx 1.000$ & $2 \shortto 3$\\
		$-$0.19225157 & $-$1.8454535 & 2.00585092 & 0.81287280 & 12.567228 & $(-)$ $\theta \approx 6.051$, $\lambda \approx 1.039$ & 3\\
		0.45463054 & $-$1.9235695 & 1.90932238 & 0.80122059 & $\approx 2 \cdot 2 \pi$ & $(-)$ $\theta \approx 5.995$, $\lambda \approx 1.045$ & 3\\
		1.51574261 & $-$2.0067614 & 1.72869916 & 0.72220662 & 12.564860 & $(-)$ $\theta \approx 5.947$, $\lambda \approx 1.052$ & 3\\
		3.26379993 & $-$2.2415051 & 1.60217526 & 0.29550907 & 12.560384 & $(-)$ $\theta \approx 6.023$, $\lambda \approx 1.000$ & $3 \shortto 2$\\
		3.40128178 & $-$2.3398344 & 1.69790628 & 0.21620397 & 12.560137 & $(-)$ $\theta_1 \approx 6.105$, $(-)$ $\theta_2 \approx 6.275$ & 2\\
		3.57582149 & $-$2.5235422 & 1.89368945 & 0 & 12.560333 & $(-)$ $\theta_p \approx 6.259$, $\lambda_s \approx 1.000$ & $2 \shortto 4$\\
		\hline $C$ & $x(0)$ & $z(0)$ & $\dot{y}(0)$ & $T/2$ & Floquet multipliers \& Krein sign & $\mu_{CZ}$\\
		\hline $-$2.78507825 & $-$2.5252333 & 0 & 3.15514897 & 12.572403 & $(-)$ $\theta_p \approx 6.259$, $\lambda_s \approx 1.000$ & $2 \shortto 4$\\
		$-$2.20230468 & $-$1.9690911 & 1.38000000 & 2.62909547 & 12.570561 & $(-)$ $\theta \approx 6.195$, $\lambda \approx 1.003$ & 3\\
		$-$0.81720281 & $-$0.7479807 & 1.73559379 & 1.56032664 & 12.568861 & $(-)$ $\theta \approx 5.864$, $\lambda \approx 1.018$ & 3\\
		0.89061923 & 0.27161684 & 1.61742872 & 0.63247393 & $\approx 2 \cdot 2 \pi$ & $(-)$ $\theta \approx 5.551$, $\lambda \approx 1.042$ & 3\\
		1.02320791 & 0.34457437 & 1.60119923 & 0.56032664 & 12.566073 & $(-)$ $\theta \approx 5.536$, $\lambda \approx 1.044$ & 3\\
		2.86941464 & 1.59425590 & 1.21232538 & $-$0.8196733 & 12.560660 & $(-)$ $\theta \approx 5.776$, $\lambda \approx 1.015$ & 3\\
		3.57582149 & 2.52462303 & 0 & $-$1.8953990 & 12.560333 & $(-)$ $\theta_p \approx 6.259$, $\lambda_s \approx 1.000$ & $2 \shortto 4$
	\end{tabular}
	\label{data_5}
\end{table}

\begin{table}[ht]\fontsize{10}{10}\selectfont \centering
	\caption{Data for two families emerging from 1:6 $\kappa_-$ (both families terminate at 1:4 $\kappa_+$).}
	\begin{tabular}{c|c|c|c|c|c|c}
		$C$ & $x(0)$ & $\dot{y}(0)$ & $\dot{z}(0)$ & $T/4$ & Floquet multipliers \& Krein sign & $\mu_{CZ}$\\
		\hline $-$3.08307598 & $-$2.9282380 & 3.51302954 & 0 & 7.856048 & $(-)$ $\theta_p \approx 6.266$, $\lambda_s \approx 1.000$ & $2 \shortto 4$\\
		$-$2.82590417 & $-$2.9274489 & 3.46849895 & 0.22200000 & 7.855771 & $(-)$ $\theta_1 \approx 6.270$, $(-)$ $\theta_2 \approx 6.274$ & 2\\
		$-$2.05755697 & $-$2.9256958 & 3.33581052 & 0.41706000 & 7.855169 & $(-)$ $\theta \approx 6.231$, $\lambda \approx 1.000$ & $2 \shortto 3$\\
		0.36834818 & $-$2.9233638 & 2.91894867 & 0.58507392 & $\approx 2 \pi$ $5/4$ & $(-)$ $\theta \approx 6.129$, $\lambda \approx 1.021$ & 3\\
		3.36677655 & $-$2.9260168 & 2.40917597 & 0.27369986 & 7.852301 & $(-)$ $\theta \approx 6.229$, $\lambda \approx 1.000$ & $3 \shortto 2$\\
		3.64903893 & $-$2.9268186 & 2.36188185 & 0.15098786 & 7.852040 & $(-)$ $\theta_1 \approx 6.265$, $(-)$ $\theta_2 \approx 6.269$ & 2\\
		3.76531791 & $-$2.9272292 & 2.34250125 & 0 & 7.851892 & $(-)$ $\theta_p \approx 6.266$, $\lambda_s \approx 1.000$ & $2 \shortto 4$\\		
		\hline $C$ & $x(0)$ & $\dot{y}(0)$ & $\dot{z}(0)$ & $T/4$ & Floquet multipliers \& Krein sign & $\mu_{CZ}$\\
		\hline $-$3.08307598 & 2.92838023 & $-$3.5132457 & 0 & 7.856048 & $(-)$ $\theta_p \approx 6.266$, $\lambda_s \approx 1.000$ & $2 \shortto 4$\\
		$-$2.86931536 & 2.92776198 & $-$3.4762696 & 0.20300000 & 7.855811 & $(-)$ $\theta \approx 6.270$, $\lambda \approx 1.006$ & 3\\
		$-$2.11373395 & 2.92627646 & $-$3.3460230 & 0.40726000 & 7.855185 & $\lambda_1 \approx 1.000$, $\lambda_2 \approx 1.049$ & $3 \shortto 4$\\
		0.31581407 & 2.92493344 & $-$2.9296442 & 0.58465849 & $\approx 2 \pi$ $5/4$ & $\lambda_1 \approx 1.021$, $\lambda_2 \approx 1.165$ & 4\\
		0.35678440 & 2.92493645 & $-$2.9226442 & 0.58466924 & 7.853965 & $\lambda_1 \approx 1.021$, $\lambda_2 \approx 1.166$ & 4\\
		3.39525998 & 2.92699047 & $-$2.4056141 & 0.26413810 & 7.852319 & $\lambda_1 \approx 1.000$, $\lambda_2 \approx 1.051$ & $4 \shortto 3$\\
		3.66107019 & 2.92741447 & $-$2.3606929 & 0.14300260 & 7.852044 & $(-)$ $\theta \approx 6.269$, $\lambda \approx 1.015$ & 3\\
		3.76531791 & 2.92763105 & $-$2.3431368 & 0 & 7.851892 & $(-)$ $\theta_p \approx 6.266$, $\lambda_s \approx 1.000$ & $2 \shortto 4$
	\end{tabular}
	\label{data_6}
\end{table}

\begin{table}[ht]\fontsize{10}{10}\selectfont \centering
	\caption{Data for branch that bifurcate from 1:1 $\kappa_+$.}
	\begin{tabular}{c|c|c|c|c|c|c}
		$C$ & $x(0)$ & $z(0)$ & $\dot{y}(0)$ & $T/2$ & Floquet multipliers \& Krein sign & $\mu_{CZ}$\\
		\hline 3.15340606 & 1.69254257 & 0 & $-$0.9511884 & 6.254196 & $\lambda_p \approx 1.473$, $\lambda_s \approx 1.000$ & $3 \shortto 4$\\
		3.10497908 & 1.64057057 & 0.32500000 & $-$0.8903832 & 6.256432 & $\lambda \approx 1.418$, $(-)$ $\theta \approx 5.997$ & 3\\
		2.50281375 & 1.14121257 & 1.02130088 & $-$0.3248755 & 6.272841 & $\lambda \approx 1.284$, $(-)$ $\theta \approx 5.405$ & 3\\
		%1.71481742 & 0.64318063 & 1.34352107 & 0.19512442 & 6.280390 & $\lambda \approx 1.207$, $(-)$ $\theta \approx 5.374$ & 3\\
		1.06209429 & 0.25735294 & 1.48471729 & 0.57281442 & $\approx 2 \pi$ & $\lambda \approx 1.157$, $(-)$ $\theta \approx 5.484$ & 3\\
		0.61412474 & $-$0.01215181 & 1.54242536 & 0.82470958 & 6.284400 & $\lambda \approx 1.126$, $(-)$ $\theta \approx 5.594$ & 3\\
		$-$0.43400738 & $-$0.81183181 & 1.69890129 & 1.69890129 & 6.285041 & $\lambda \approx 1.193$, $(-)$ $\theta \approx 6.282$ & 3\\
		\hdashline
		$-$0.43400072 & $-$0.81183281 & 1.69891361 & 1.46801799 & 6.285041 & $\lambda_1 \approx 1.193$, $\lambda_2 \approx 1.001$ & 4\\
		$-$0.33553943 & $-$0.80429683 & 1.84771240 & 1.40521652 & $\approx 2 \pi$ & $\lambda_1 \approx 1.516$, $\lambda_2 \approx 1.038$ & 4\\
		0.73611388 & 0.11872418 & 2.40681480 & 0.32732843 & 6.274234 & $\lambda_1 \approx 2.749$, $\lambda_2 \approx 1.044$ & 4\\
		1.61080283 & 1.64821787 & 2.23796057 & $-$1.3510095 & 6.264739 & $\lambda_1 \approx 3.760$, $\lambda_2 \approx 1.001$ & 4\\
		\hdashline 1.61080283 & 1.64821872 & 2.23796038 & $-$1.3510105 & 6.264739 & $\lambda \approx 3.760$, $(-)$ $\theta \approx 6.283$ & 3\\
		1.39262658 & 1.73365484 & 2.37502843 & $-$1.5142517 & 6.258164 & $\lambda \approx 3.376$, $(-)$ $\theta \approx 6.243$ & 3\\
		0.73257748 & 0.67443663 & 3.04798979 & $-$1.6021384 & 6.211295 & $\lambda \approx 2.178$, $(-)$ $\theta \approx 6.122$ & 3\\
		0.66946846 & 0.37443663 & 3.09929510 & $-$1.3331532 & 6.159828 & $\lambda \approx 1.705$, $(-)$ $\theta \approx 6.025$ & 3\\
	\end{tabular}
	\label{data_7}
\end{table}

\end{appendix}

\clearpage

\noindent
\textbf{Acknowledgement.}\ 
The author acknowledges support by the Deutsche Forschungsgemeinschaft (DFG, German Research Foundation), Project-ID 541062288.\ He also wishes to thank Urs Frauenfelder for helpful discussions.\ He also would like to thank Alexander Batkhin and Alexander Bruno for providing English annotations from some of Bruno's preprints in Russian.\ Furthermore, he gratefully acknowledges Otto van Koert for valuable discussions, his warm hospitality and fostering an inspiring atmosphere during the author's stay as a visiting researcher at the Seoul National University in November 2025, where significant parts of the section on the analytical approach to comet orbits were developed.\

%\clearpage

\end{document}